\documentclass[reqno,12pt]{article}
\usepackage{amssymb,amsthm,amsmath,amsfonts}
\usepackage{epsfig}

\theoremstyle{plain}
\begingroup

\newtheorem{thm}{Theorem}[section]
\newtheorem{theorem}{Theorem}[section]
\newtheorem{cor}[thm]{Corollary}
\newtheorem{corollary}[thm]{Corollary}
\newtheorem{lemma}[thm]{Lemma}

\newtheorem{proposition}[thm]{Proposition}

\endgroup
\theoremstyle{definition}
\newtheorem{defn}{Definition}[section]

\newtheorem{remark}[defn]{Remark}

\newtheorem{example}[defn]{Example}

\theoremstyle{remark}

\numberwithin{equation}{section}
\numberwithin{figure}{section}

\DeclareMathOperator{\re}{Re} \DeclareMathOperator{\im}{Im}

\DeclareMathOperator*{\res}{\mathrm{Res}}

\def\I{\mathrm{i}}

\def\D{{\mathbb D}}

\def\C{{\mathbb C}}
\def\P{{\mathbb P}}

\begin{document}

\title{
On the dynamics of roots and poles for solutions of the
Polubarinova-Galin equation }
\author{
Bj\"orn Gustafsson\textsuperscript{1}
Yu-Lin Lin\textsuperscript{2}\\
}

\date{December 1, 2011}

\maketitle

\begin{abstract}
We study the dynamics of roots of $f'(\zeta,t)$, where $f(\zeta,t)$
is a locally univalent polynomial solution of the Polubarinova-Galin equation
for the evolution of the conformal map onto a Hele-Shaw blob subject to injection
at one point. We give examples of the sometimes complicated motion of roots, but show also that
the asymptotic behavior is simple.
More generally we allow $f'(\zeta,t)$ to be a rational function and give sharp estimates for the
motion of poles and for the decay of the Taylor coefficients. We also prove that any global in time
locally univalent solution actually has to be univalent.
\end{abstract}

 \footnotetext[1]
{Department of Mathematics, KTH, 100 44, Stockholm, Sweden.\\
Email: \tt{gbjorn@kth.se}}
\footnotetext[2]
{Department of Mathematics, KTH, 100 44, Stockholm, Sweden.\\
Email: \tt{ylli@kth.se}}

\noindent {\it Acknowledgements:} This work has been supported by visiting grants from
the Royal Institute of Technology (Stockholm) and Academia Sinica at Taipei, and in addition by
the European Science Foundation Research Networking Programme HCAA, the Mittag-Leffler Institute
(Djursholm, Sweden) and the G\"oran Gustafsson Foundation,

The authors are grateful to Saleh Tanveer for illuminating discussions
and to Govind Menon, who helped the second author start out on the subject.


\section{Introduction}

Polynomial, rational and logarithmic solutions of the Polubarinova-Galin equation for
the conformal map onto a growing or shrinking Hele-Shaw blob of viscous fluid
have been studied in many papers the last few decades, see for example
\cite{Shraiman-Bensimon84},
\cite{Gustafsson84},
\cite{Howison86},
\cite{Mineev90},
\cite{Tanveer93},
\cite{Abanov-etal09},
\cite{Hohlov-Howison93},
\cite{Dawson-Mineev94},
\cite{Huntingford95},
\cite{Richardson97},
\cite{Vondenhoff08},
\cite{Lin09a},
\cite{Lin09b},
and also \cite{Gustafsson-Vasiliev06}.
In \cite{Shraiman-Bensimon84} D.~Bensimon and B.~Shraiman set up the dynamical equations,
in the polynomial case, for how the zeros of the derivative of the conformal map move
in the complex plane, and they proposed studying the dynamics of the zeros.
As far as we know nobody has yet undertaken this task to any substantial extent.
The purpose of the present paper is to start investigations in such a direction.

The general setting is that $f(\cdot,t):\D\to \Omega(t)\subset \C$ is the time dependent
normalized ($f(0,t)=0$, $f'(0,t)>0$) conformal map onto the fluid region $\Omega(t)$, which
has a source of strength $q(t)>0$ at the origin. The evolution is then described by
the Polubarinova-Galin equation
\begin{equation}\label{pg0}
{\rm Re}\left[\dot{f}(\zeta,t)\overline{\zeta
f'(\zeta,t)}\right]=q(t) \quad {\rm for}\,\,\zeta\in\partial
\mathbb{D}.
\end{equation}
We focus mainly on the case that the derivative $g(\zeta,t)=f'(\zeta,t)$ is a rational function,
a property which is preserved in time under (\ref{pg0}). Then $g(\zeta,t)$, and hence $f(\zeta,t)$,
is completely determined
by the zeros $\omega_1(t), \dots, \omega_m(t)$ and poles $\zeta_1(t),\dots, \zeta_n(t)$ of $g(\zeta,t)$,
both located outside the unit disk, together with an overall scale factor $b(t)$:
\begin{equation}\label{structureg0}
g(\zeta,t)
=b(t)\frac{\prod_{k=1}^{m}(\zeta-\omega_{k}(t))}{\prod_{j=1}^{n}(\zeta-\zeta_{j}(t))}.
\end{equation}

The motion of the poles $\zeta_j(t)$ has previously been studied, for example it is known \cite{Tanveer93},
\cite{Gustafsson-Vasiliev06} that they
always move away from the origin. In the present paper we obtain more precise
estimates of their speed and locations.
The zeros $\omega_k(t)$ show up a more complicated behavior than the poles, for example they may occasionally move towards the origin, and even reach the unit circle, with possible break down of the solution as the result. Zeros may also collide with each other or escape to infinity in finite time.
However the asymptotic behavior is simple: if the solution does not break down, then as $t\to \infty$ the zeros move to infinity arranged asymptotically as the corners of a regular polygon of growing size.

Besides the zeros and poles we also study the asymptotics of the coefficients in the power series
$$
f(\zeta,t)= \sum_{j=1}^\infty a_j(t)\zeta^j
$$
in the case that $g(\zeta,t)$ is rational. The leading coefficient $a_1(t)>0$ increases to infinity as $t\to \infty$,
in fact can be normalized for example so that $a_1(t)=e^t$ (which then fixes $q(t)$), while the others
tend to zero very quickly:
\begin{equation}\label{decay0}
a_j (t)=O(\frac{1}{a_1(t)^j}) \quad (j\geq 2)
\end{equation}
as $t\to\infty$.

Much of the behavior of solutions of the Polubarinova-Galin equation (\ref{pg0}) can be understood in terms
of harmonic moments and quadrature identities. These provide enough conserved quantities
to make the Hele-Shaw problem fully integrable, hence in principle
algebraically solvable. In practice it is not that easy because of the degree of complexity of
the integrals of motion when expressed for example in terms of the Taylor coefficients of $f(\zeta,t)$.
However the conservation laws allow for natural concepts of weak solutions, and with such relaxed forms of
solutions it is possible to let zeros of $g(\zeta,t)$ penetrate the unit circle and enter the unit
disk. Non locally univalent solutions of this kind will be studied in a forthcoming paper \cite{Gustafsson-Lin11},
and part of the aim of the present paper is to set the stage for these future investigations.

The organization of the paper is as follows. Section~\ref{sec:prep} contains all necessary preliminaries, in particular we set up the notations to be used. In Section~\ref{sec:dynamics} we derive the dynamical equations for the zeros and poles of $g(\zeta,t)$ by first writing (\ref{pg0}) as an equation in $g$ alone (see (\ref{dotlogg})) and then identifying the residues
in this formula (Theorem~\ref{dynamics}). When everything is spelled out (Theorem~\ref{thm:dynamics} and
\ref{thm:dynamics1}) one gets a rather involved system of ordinary
differential equations for $\omega_k(t)$, $\zeta_j(t)$, $b(t)$. From Theorem~\ref{thm:dynamics1} we obtain precise estimates of the speeds of the poles
(Corollary~\ref{cor:polesmoveout}), and a kind of conservation law for the dynamics (Proposition~\ref{prop:conservation}).

The Polubariova-Galin equation (\ref{pg0}) is not sensitive for loss of univalence of the solution
as long as it stays locally univalent. However, loss of univalence is always followed by
loss also of local univalence and break down of the solution at a later, but finite, time
(Theorem~\ref{lemma2}).
In \cite{Gustafsson-Lin11} we will show how it in such cases still can be continued
as a kind of weak solution spreading on a branched Riemann surface.

In Section~\ref{sec:asymptotic} we give several examples (Example~\ref{ex:collision}, \ref{ex:huntingford},
\ref{ex:onecloseroot}) of the behavior of roots $\omega_k(t)$ in the polynomial case, for example collision of roots,
plus two theorems asserting the previously mentioned asymptotic behavior.
In the first (Theorem~\ref{thmgg}) it is  assumed {\it a priori} that the solution is global in time, while in the second
(Theorem~\ref{tem}) only sufficiently strong assumptions on the initial data are made.
Section~\ref{sec:multiplecut}, finally, contains precise estimates for the poles $\zeta_j(t)$ in the rational case
(Theorem~\ref{first}) as well as the asymptotic estimates (\ref{decay0}) of the Taylor coefficients of $f(\zeta,t)$
(Theorem~\ref{third}; see also Lemma~\ref{decay}). There are also some sharper estimates (e.g.,
Corollary~\ref{four}) for the case that the sequence of harmonic moments contains gaps.


\section{Preparatory material}\label{sec:prep}

\subsection{List of notations}\label{sec:notations}

We first list some general notations which will be used in the
paper.
\begin{itemize}

\item $\mathbb{D}=\{\zeta\in \mathbb{C}:|\zeta|< 1\}$,
$\mathbb{D}(a,r)=\{\zeta\in \mathbb{C}:|\zeta-a|< r\}$.

\item $dm= dm(z)=dxdy=\frac{1}{2\I}d\bar{z}dz$
($z=x+\I y$), area measure in the $z$ plane.

\item $\omega^{*}=\frac{1}{\overline{\omega}}$, for $\omega\in \mathbb{C}$.

\item $h^{*}(\zeta)=\overline{h(1/\bar{\zeta})}
=\sum_{j=1}^{m}\overline{b}_{j}\zeta^{-j}$, where
$h(\zeta,t)=\sum_{j=1}^{m}b_{j}\zeta^{j}$.

(There is a slight ambiguity in this notation: we have $\zeta^*=1/\bar\zeta$
if $\zeta$ is considered as a point, whereas $f^*(\zeta)=1/\zeta$ for the
function $f(\zeta)=\zeta$.)

\item $\dot{f}(\zeta,t)=\frac{\partial}{\partial t}f(\zeta,t)$,
$f^{'}(\zeta,t)=\frac{\partial}{\partial\zeta}f(\zeta,t)$.

\item With ${E}\subset\mathbb{C}$ any set which contains the origin,
\begin{align}
\mathcal{O}({E})=&\{f:  \mbox{$f$ is analytic in some neighborhood
of
$E$}\},\notag\\
\mathcal{O}_{\rm norm}({E})=&\{f\in \mathcal{O}({E}):  f(0)=0,
f'(0)>0\},\notag\\
\mathcal{O}_{\rm locu}({E})=&\{f\in \mathcal{O}_{\rm norm}({E}):
f'\neq
0 \,\,\mbox{on}\,\,E\},\notag\\
\mathcal{O}_{\rm univ}({E})=&\{f\in \mathcal{O}_{\rm locu}({E}):
\mbox{$f$ is univalent (one-to-one) on ${E}$}\}.\notag
\end{align}

\item If $a(t)$, $b(t)$ are positive functions
$a(t)\sim b(t)$ will mean that there exist constants $0<c<C<\infty$
such that
$$
c\leq \frac{a(t)}{b(t)}\leq C
$$
for all $t$.

\end{itemize}\par


\subsection{Basic equations}\label{sec:basicequations}

In the paper we shall study certain aspects of the {\bf
Polubarinova-Galin equation}, which describes the evolution of a
simply connected Hele-Shaw fluid blob in the complex plane driven by
injection or suction at one point, chosen to be the origin. A smooth
map $t\mapsto f(\cdot,t)\in \mathcal{O}_{\rm
univ}(\overline{\mathbb{D}})$ will be called a \textbf{(univalent)
solution} of the Polubarinova-Galin
equation if it satisfies
\begin{equation}\label{pg1}
{\rm Re}\left[\dot{f}(\zeta,t)\overline{\zeta
f'(\zeta,t)}\right]=q(t) \quad {\rm for}\,\,\zeta\in\partial
\mathbb{D}
\end{equation}
in the pointwise sense.
Here $q(t)$ is a real-valued continuous function, which is given in
advance and which represents the strength of the source/sink at the
origin. Typically $q=\pm 1$, which corresponds to injection (plus
sign) or suction (minus sign) at a rate $2\pi$. Since the
transformation $t\mapsto -t$ changes $q$ to $-q$ in (\ref{pg1}) it
is enough to discuss one of the cases $q>0$ and $q<0$. In general we
shall take $q>0$.

Equation (\ref{pg1}) expresses that the image domains
$\Omega(t)=f(\mathbb{D},t)$ evolve in such a way that
\begin{equation}\label{lg}
\frac{d}{dt}\int_{\Omega(t)} hdm =2\pi{q(t)}h(0)
\end{equation}
for every function $h$ harmonic in a neighborhood of
$\overline{\Omega(t)}$. This means that the speed of the boundary
$\partial\Omega(t)$ in the normal direction equals $q(t)$ times the
normal derivative of the Green's function of $\Omega(t)$ with a pole
at $z=0$. The terminology {\it Laplacian growth} is also used for
this kind of evolution. The equivalence between (\ref{pg1}) and
(\ref{lg}) follows from the general formula
\begin{equation}\label{generalevolution}
\frac{d}{dt}\int_{\Omega(t)} \varphi dm
=\int_{\partial\D}\varphi(f(\zeta,t)){\rm
Re}\left[\dot{f}(\zeta,t)\overline{\zeta f'(\zeta,t)}\right]d\theta
\quad (\zeta=e^{\I\theta}),
\end{equation}
valid for any smooth evolution $t\mapsto f(\cdot,t)\in
\mathcal{O}_{\rm univ}(\overline{\mathbb{D}})$ and for any smooth
test function $\varphi$ in the complex plane. The derivation of (\ref{generalevolution})
is straightforward (and omitted).

On choosing $h(z)=z^k$, $k=0,1,2,\dots$ in (\ref{lg}) it follows
that the harmonic moments
\begin{equation}\label{defmoments}
M_{k}(t)=\frac{1}{\pi}\int_{\Omega(t)}z^{k}dm = \frac{1}{2\pi
\I}\int_{\partial\mathbb{D}}f(\zeta,t)^k f^*(\zeta,t) f'(\zeta,t)d\zeta
\end{equation}
are conserved quantities, except for the first one,
which is related to $q$ by $\frac{d}{dt}{M}_0(t)=2q(t)$. Thus
$$
M_0(t)=M_0(0)+2Q(t),
$$
where $Q(t)$ is the accumulated source up to time $t>0$:
\begin{equation}\label{Q}
Q(t)=\int_0^t q(s)ds.
\end{equation}
The preservation of $M_1, M_2, \dots$ characterizes Laplacian growth
in the simply connected case.

One may consider the equation (\ref{pg1}) on different levels of
generality. It is natural to keep the normalization $f(0)=0$,
$f'(0)>0$, in fact the coupling to (\ref{lg}) depends on this, but
(\ref{pg1}) makes sense for any $f\in\mathcal{O}_{\rm
norm}(\overline{\mathbb{D}})$, at least as long one makes sure that
$q(t)=0$ whenever a zero of $f'$ appears on $\partial\mathbb{D}$. In
a forthcoming paper \cite{Gustafsson-Lin11} we shall deal with this general case,
while in the present paper we shall only consider locally univalent
functions, $f\in\mathcal{O}_{\rm locu}(\overline{\mathbb{D}})$. The
mathematical treatment of (\ref{pg1}) in this case is exactly the
same as in the `physical' case $f\in\mathcal{O}_{\rm
univ}(\overline{\mathbb{D}})$. We shall then speak of a
\textbf{locally univalent solution} of the Polubarinova-Galin
equation.

When $f\in\mathcal{O}_{\rm locu}(\overline{\mathbb{D}})$, then
$\dot{f}/\zeta f^{'}\in \mathcal{O}(\overline{\mathbb{D}})$ and
equation (\ref{pg1}) can be solved for $\dot{f}$ by dividing both
sides with $|\zeta f'|^2$. The result is an equation which we shall refer
to as the {\bf L\"owner-Kufarev equation}, namely
\begin{equation}\label{lk}
\dot{f}(\zeta,t)=\zeta f'(\zeta,t)P(\zeta,t)
\quad(\zeta\in\mathbb{D}),
\end{equation}
where $P(\zeta,t)$ is the analytic function in $\D$ whose real part
has boundary value $q(t)|f'(\zeta,t)|^{-2}$ and which is normalized
by $\im P(0,t)=0$. Explicitly $P(\zeta,t)$ is given by
\begin{equation}\label{poisson}
P(\zeta,t) =\frac{1}{2\pi i}
\int_{\partial\D}\frac{q(t)}{|f'(z,t)|^2}\frac{z+\zeta}{z-\zeta}\frac{dz}{z}
\quad(\zeta\in\mathbb{D}).
\end{equation}
The right member in (\ref{poisson}) is the Poisson-Schwarz integral
\begin{equation}\label{poisson-schwarz}
\mathcal{P}_\mu (\zeta)=\int_0^{2\pi}
\frac{e^{i\theta}+\zeta}{e^{i\theta}-\zeta}d\mu(\theta)
\quad(\zeta\in\mathbb{D})
\end{equation}
for the measure
\begin{equation}\label{measure}
d\mu(\theta,t)=
\frac{q(t)}{|f'(e^{i\theta},t)|^{2}}\frac{d\theta}{2\pi}.
\end{equation}
If $q(t)>0$ then $\mu$ is positive and $\re P>0$ in $\D$.

\begin{remark} In the general literature on univalent functions,
``L\"owner-Kufarev equation'' usually refers to any equation of the
type (\ref{lk}) without any particular coupling between $f$ and $P$.
The Hele-Shaw case is characterized by the feed-back relation
(\ref{poisson}) between  $f$ and $P$.
\end{remark}

Expanding $f$ in a power series,
\begin{equation}\label{powerseries}
f(\zeta,t)= a_1(t)\zeta +a_2(t)\zeta^2+\dots,
\end{equation}
it follows from (\ref{lk}) that
\begin{equation}\label{dotaovera}
\frac{\dot{a}_1(t)}{a_1(t)}=P(0,t)=\int d\mu(\cdot,t).
\end{equation}
In particular, if $q(t)$ is chosen so that $\mu$ becomes a
probability measure, namely
\begin{equation}\label{muprobability}
\frac{1}{q(t)}=\frac{1}{2\pi}\int_0^{2\pi}\frac{d\theta}{|f'(e^{\I\theta},t)|^{2}},
\end{equation}
then $P(0,t)=1$ and $a_1(t)=a_1(0)e^t$.

When $f\in\mathcal{O}_{\rm locu}(\overline{\mathbb{D}})$ then
$P\in\mathcal{O}(\overline{\mathbb{D}})$, in fact the right member
of (\ref{lk}) extends analytically as far as $f$ does (see
\cite{Gustafsson84}). We shall keep the notation $P=P(\zeta,t)$ also for
the analytic extension of the Poisson integral beyond
$\overline{\mathbb{D}}$.

As a general notation throughout the paper, we set
\begin{equation}\label{derivative}
g(\zeta,t)=f' (\zeta,t).
\end{equation}
The function $g$ in fact turns out to be more fundamental than $f$
itself. Of course, $f$ can be recaptured from $g$ by
$$
f(z,t)=\int_0^z g(\zeta,t)d\zeta.
$$

Most of the paper will deal with the case that $g$ is a rational
function, or perhaps better to say, $g(\zeta)d\zeta$ is a rational
differential, in other words an Abelian differential on the Riemann
sphere. If $g$ has residues then $f$ will have logarithmic poles,
besides ordinary poles. The terminology {\bf Abelian domain} for the image
domain $\Omega=f(\mathbb{D})$ has been used \cite{Varchenko-Etingof92}
for this case. Alternatively one may speak of $\Omega$ being a
{\bf quadrature domain} (see \cite{Gustafsson-Shapiro05} for
the terminology), which in the present case means that a finite quadrature identity of
the kind
\begin{equation}\label{qi}
\int_\Omega h(z)dxdy=\sum_{j=1}^r c_j\int_{\gamma_j}h(z)dz
+\sum_{j=0}^\ell\sum_{k=1}^{n_j-1}a_{jk}h^{(k-1)}(z_j)
\end{equation}
holds for integrable analytic functions $h$ in $\Omega$. Here the points
$z_j\in\Omega$ are fixed (i.e., independent of $h$), with specifically $z_0=0$,
the $c_j$, $a_{jk}$ are fixed coefficients, and the
$\gamma_j$ are arcs in $\Omega$ with end points among the $z_j$.
This sort of structure is
stable under Hele-Shaw flow because, as is seen from (\ref{lg}),
what happens under the evolution is only that the right member is
augmented by the term $2\pi Q(t) h(0)$.

When $g=f'$ is rational we shall write it on the form
\begin{equation}\label{structureg}
g(\zeta,t)
=b(t)\frac{\prod_{k=1}^{m}(\zeta-\omega_{k}(t))}{\prod_{j=1}^{n}(\zeta-\zeta_{j}(t))}
=b(t)\frac{\prod_{i=1}^{m}(\zeta-\omega_{i}(t))}{\prod_{j=1}^{\ell}(\zeta-\zeta_{j}(t))^{n_j}}.
\end{equation}
Here $m\geq n=\sum_{j=1}^\ell {n_j}$, $|\zeta_j|>1$ and repetitions
are allowed among the $\omega_k$, $\zeta_j$ to account for multiple
zeros and poles. Then, with the argument of $b(t)$ chosen so that
$g(0,t)>0$, $f\in\mathcal{O}_{\rm locu}(\overline{\mathbb{D}})$ if
and only if $|\omega_k|>1$, $|\zeta_j|>1$ for all $k$ and $j$. The
assumption $m\geq n$ means that $g(\zeta)d\zeta$, as a differential,
has at least a double pole at infinity, which the Hele-Shaw
evolution in any case will force it to have because the source/sink
at the origin creates a pole of $f$ at infinity.

The form (\ref{structureg}) is stable in time, with the sole exception that
when $m=n$ the pole of $f$ may disappear at one moment of time. Then the value
of $m$ drops below $n$ at this moment. This can easily be explained in terms of
the quadrature identity (\ref{qi}): under the Hele-Shaw evolution there will be
one time dependent term, namely $(a_{01}+2\pi Q(t))h(0)$, and this may vanish for
one value of $t$. If there are no terms $a_{0k}h^{(k-1)}(0)$ with $k\geq 2$
this will cause $f$ to temporarily lose its pole at infinity. The phenomenon
is illustrated in Example~\ref{ex:drop} below.

The rightmost member of (\ref{structureg}) will be used when we need
to be explicit about the orders of the poles. The convention then is
that $\zeta_1,\dots,\zeta_\ell$ are distinct and $n_j\geq 1$. Thus
$n=\sum_{j=1}^\ell {n_j}$, and in the full sequence
$\zeta_1,\dots,\zeta_n$, the tail $\zeta_{\ell+1},\dots,\zeta_n$
will be repetitions of (some of) the $\zeta_1,\dots,\zeta_\ell$
according to their orders. In equations (\ref{qi}) and
(\ref{structureg}), $\ell$ and the $n_j$ are the same.

For later use we here also exhibit $\log g$ and its derivatives:
\begin{equation}\label{logg}
\log g(\zeta,t) =\log
b(t)+{\sum_{k=1}^{m}\log({\zeta-\omega_{k}(t)}})
-{\sum_{j=1}^{n}\log({\zeta-\zeta_{j}(t)}}),
\end{equation}
\begin{equation}\label{loggdot}
\frac{\partial}{\partial t}\log g(\zeta,t)=\frac{\dot{b}(t)}{b(t)}
-{\sum_{k=1}^{m}\frac{\dot{\omega}_k(t)}{\zeta-\omega_{k}(t)}}
+{\sum_{j=1}^{n}\frac{\dot{\zeta}_j(t)}{\zeta-\zeta_{j}(t)}},
\end{equation}
\begin{equation}\label{loggprime}
\frac{\partial}{\partial \zeta}\log g(\zeta,t)=
{\sum_{k=1}^{m}\frac{1}{\zeta-\omega_{k}(t)}}
-{\sum_{j=1}^{n}\frac{1}{\zeta-\zeta_{j}(t)}}.
\end{equation}
We remark that $\log|g|$ can be viewed as the logarithmic potential
of the charge distribution which puts positive unit charges at the
zeros of $g$, negative unit charges at the poles (added up according
to multiplicities). This is the charge distribution corresponding to
the divisor $(g)$ of $g$, which is defined as the formal linear
combination
\begin{equation}\label{divisor}
(g)=\sum_{k=1}^m 1\cdot (\omega_k)-\sum_{j=1}^n 1\cdot
(\zeta_j)-(m-n)\cdot (\infty).
\end{equation}
(This is the divisor of $g$ as a function, the divisor of $gd\zeta$
as a differential has two additional poles at infinity.) One main
subject of the present paper is the study of the dynamics of $(g)$
under Hele-Shaw flow.

Note that $\log|g|$, and even its restriction to $\partial\D$,
contains all information of $f$. In fact, given $u=\log|g|$ on
$\partial \D$ we can extend it harmonically to $\D$, then form its
harmonic conjugate $v$, normalized by $v(0)=0$, and finally define
$f$ by
$$
f(\zeta)=\int_0^\zeta \exp({u(z)+\I v(z)})dz.
$$
Note also that $(\log|g|)|_{\partial\D}$ is a ``free'' function,
i.e., is subject to no constraints besides regularity (real
analyticity is needed to start a Hele-Shaw evolution in a full
neighborhood of $t=0$, see e.g. \cite{Reissig93}).

The structure (\ref{structureg}) of $g$ means that $f$ is of the
form
\begin{equation}\label{structuref}
f(\zeta,t)=\sum_{j=1}^\ell e_{j}\log(\zeta-\zeta_j(t))
+\sum_{j=1}^\ell\sum_{k=1}^{n_j-1}\frac{c_{jk}(t)}{(\zeta-\zeta_j(t))^k}
+\sum_{k=0}^{m-n+1} d_k(t) \zeta^k.
\end{equation}
The coefficients $e_{j}$, which are the residues of $g(\zeta)d\zeta$, will
not depend on time (this is a consequence of (\ref{pg1})). If they
are not all zero, then $f$ is single-valued only outside a suitable
system of `cuts' in $\C\setminus \overline{\D}$ connecting the
logarithmic poles. Note that if $\sum_{j=1}^\ell e_{j}\ne 0$ then
$f$ automatically has a logarithmic pole at infinity, hence one of
the cuts has to reach infinity in this case. The relationship
between (\ref{structuref}) and (\ref{qi}) for $\Omega=\Omega
(t)=f(\D,t)$ is that $z_j=f(\zeta_j^*(t),t)$, $1\leq j\leq \ell$,
$z_0=f(\infty^*,t)=f(0,t)=0$, and
that the above mentioned `cuts' can be taken to be the reflections
in the unit circle of the arcs $f^{-1}(\gamma_j)$. Moreover,
$\sum_{j=1}^\ell n_j=n$ and $n_0=m-n+2$, hence $\sum_{j=0}^\ell n_j=m+2$.
One may also note that $b=(m-n+1)d_{m-n+1}$.
The two systems of coefficients, $\{c_j, a_{jk}\}$ and $\{e_j, c_{jk}, d_k\}$,
are related by nonlinear equations (see e.g.
\cite{Varchenko-Etingof92}).


\section{Dynamics of zeros and poles in the rational case}
\label{sec:dynamics}

In this section we set up the dynamical equations for zeros and
poles in the locally univalent rational case. We take $f$ and $g$ to
be of the form (\ref{structuref}), (\ref{structureg}), so that the
divisor of $g$ is given by (\ref{divisor}). When $g$ is rational
also $P$ is rational. When we refer to $P$ outside $\D$ we shall
always mean this rational function, or more generally the analytic
continuation of $P|_\D$ (when it exists). First we write the
dynamical equations in a general form.

\begin{theorem}
Under the above assumptions, the evolution of the divisor $(g)$
under Hele-Shaw flow governed by the L\"owner-Kufarev equation
(\ref{lk}) is given by
\begin{align}\label{dynamics}
\dot{\omega}_k(t) &=-{\res}_{\zeta=\omega_k} \,
\left[\zeta P(\zeta,t)\frac{g'(\zeta,t)}{g(\zeta,t)}\right],\\
\dot{\zeta}_j(t) &={\res}_{\zeta=\zeta_j} \,\left[ \zeta
P(\zeta,t)\frac{g'(\zeta,t)}{g(\zeta,t)}\right],
\end{align}
valid whenever $\omega_k$, $\zeta_j$ is a simple zero/pole
(respectively) of $g$.

If $\omega_k$ is a multiple zero, with say $\omega_\ell=\omega_k$
for $\ell$ in an index set $K\subset\{1,2,,\dots, m\}$ (containing
$k$), then the individual $\dot{\omega}_k(t)$ need not exist, but we
have instead
$$
\frac{d}{dt}\sum_{\ell\in K}{\omega}_\ell (t)
=-{\res}_{\zeta=\omega_k} \, \left[\zeta
P(\zeta,t)\frac{g^{'}(\zeta,t)}{g(\zeta,t)}\right].
$$
For multiple poles the $\dot{\zeta}_j$ do exist, and we have
$$
n_j\dot{\zeta}_j(t) ={\res}_{\zeta=\zeta_j} \,\left[ \zeta
P(\zeta,t)\frac{g'(\zeta,t)}{g(\zeta,t)}\right]
$$
for $1\leq j\leq \ell$.
\end{theorem}

\begin{proof}
From (\ref{lk}), i.e., $\dot{f}=\zeta g P$, we get $\dot{g}=(\zeta g
P)^{'}$ and hence
\begin{equation}\label{dotlogg}
\frac{\partial}{\partial t}\log g (\zeta,t)= \zeta
P(\zeta,t)\frac{\partial}{\partial \zeta}\log g (\zeta,t)  +
\frac{\partial}{\partial \zeta}(\zeta P(\zeta,t)).
\end{equation}
At first this equality holds in $\D$, but since both members are
rational functions it becomes an identity between two rational
functions.

From (\ref{loggdot}) we see that the $\dot{\omega}_k$ and
$\dot{\zeta}_j$ are the residues of $\frac{\partial}{\partial t}\log
g$. Since the last term in (\ref{dotlogg}), being a total
derivative, has no residues we immediately get the equations for
$\dot{\omega}_k$ and $\dot{\zeta}_j$. The statement concerning
multiple roots and poles also follows easily.
\end{proof}

\begin{remark}
It should be mentioned that multiple zeros never survive for any
period of time, only collisions can occur. Multiple poles are
however stable, they never split or collide.
\end{remark}

In addition to depending on the zeros and poles, $g$ in
(\ref{structureg}) also depends on the factor $b$. The $m+n+1$
complex parameters $\omega_1,\dots,\omega_m$,
$\zeta_1,\dots,\zeta_n$, $b$ are connected to the coefficient
$a_1(t)=g(0,t)>0$ by
\begin{equation}\label{parameters}
(-1)^{m-n}b\,\frac{\prod_{k=1}^m \omega_k}{\prod_{j=1}^n
\zeta_j}=a_1.
\end{equation}
In particular, the imaginary part of the left member vanishes, which
means that the mentioned parameters are subject to one real
constraint. Taking the logarithmic time derivative of
(\ref{parameters}) and using (\ref{dotaovera}) give an evolution
equation for $b(t)$:
\begin{equation}\label{conservation}
\frac{d}{dt}\log b(t)= \int d\mu (t)-\frac{d}{dt}\sum_{k=1}^m\log
\omega_k(t)+\frac{d}{dt}\sum_{j=1}^n\log \zeta_j(t).
\end{equation}

We now proceed to evaluate the Poisson integral $P(\zeta,t)$ in
(\ref{dynamics}). This can be done by a simple residue calculus in
(\ref{poisson}), using that $|g(\zeta,t)|^2=g(\zeta,t)g^*(\zeta,t)$
is a rational function in $\zeta$. The calculation becomes even more
transparent if everything is done at an algebraic level, by which it
essentially reduces to an expansion in partial fractions.

The rational function $q(t)/g(\zeta,t)g^*(\zeta,t)$ has
poles at the zeros of $g$ and $g^*$, i.e., at $\omega_1,
\dots,\omega_m, \omega_1^*, \dots,\omega_m^*$. At infinity it has
the behavior (by (\ref{structureg}))
\begin{align}\label{C}
\lim_{\zeta\to\infty}\frac{q(t)}{g(\zeta,t)g^*(\zeta,t)}
=A_\infty=\begin{cases}
\frac{q\prod_{j=1}^n\bar{\zeta}_j}{|b|^2\prod_{j=1}^m\bar{\omega}_j}=\frac{q}{ba_1}
\quad &{\rm if\,\,} m=n, \\
0 \quad &{\rm if\,\,} m>n.
\end{cases}
\end{align}
Assuming for simplicity that the roots $\omega_1, \dots,\omega_m$
are distinct it follows that
\begin{equation}\label{partialfractions}
\frac{q(t)}{g(\zeta,t)g^*(\zeta,t)} =A_\infty+\sum_{k=1}^m
\frac{\overline{A}_k}{\overline{\omega}_k}+\sum_{k=1}^m
\left[\frac{A_k}{\zeta-\omega_k}
+{\frac{\overline{A}_k\zeta}{1-\overline{\omega}_k\zeta}}\right],
\end{equation}
where the coefficients
$A_k=A_k(q,b,\omega_1,\dots,\omega_m,\zeta_1,\dots,\zeta_n)$ are
given by
\begin{equation}\label{Ak}
A_k=\frac{q(t)}{g'(\omega_k,t)g^*(\omega_k,t)}
=\frac{q}{|b|^2}\cdot\frac{\prod_{j}(\omega_{k}-\zeta_{j})\prod_{j}\overline{(\omega_{k}^{*}-\zeta_{j})}}
{\prod_{j\neq
k}(\omega_{k}-\omega_{j})\prod_{j}\overline{(\omega_{k}^{*}-\omega_{j})}}
\end{equation}
for $1\leq k\leq m$.
Note that $A_k\ne 0$ (when $f\in\mathcal{O}_{\rm locu}(\overline{\D})$).

Now, $P(\zeta,t)$ is by definition (\ref{poisson}) that holomorphic
function in $\mathbb{D}$ whose real part has boundary values
$q(t)/g(\zeta,t)g^*(\zeta,t)$ and whose imaginary part vanishes at
the origin. The function (\ref{partialfractions}) itself certainly
has the right boundary behaviour on $\partial\mathbb{D}$, but it is
not holomorphic in $\mathbb{D}$. On the other hand, the two types of
polar parts occurring in (\ref{partialfractions}) have the same real
parts on the boundary:
$$
{\rm Re\,}\frac{A_k}{\zeta-\omega_k}={\rm
Re\,}{\frac{\overline{A}_k\zeta}{1-\overline{\omega}_k\zeta}} \quad
{\rm on\,\,}\partial\mathbb{D}.
$$
Therefore, without changing the real part on the boundary we can
make the function (\ref{partialfractions}) holomorphic in
$\mathbb{D}$ by a simple exchange of polar parts. In addition, one
can add a purely imaginary constant to account for the normalization
of $P$ at the origin. The result is that
\begin{equation}\label{P}
P(\zeta,t)=A_0+\sum_{j=1}^m \frac{2A_j}{\zeta-\omega_j}
\end{equation}
for a suitable constant $A_0$. Since $P(0,t)=\int d\mu$ is real, the
imaginary part of $A_0$ is given by
\begin{equation}\label{imA0}
\im A_0=\im \sum_{j=1}^m \frac{2A_j}{\omega_j}.
\end{equation}

For the reflected kernel
$P^*(\zeta,t)=\overline{P({{1}/{\bar\zeta}},t)}$ we get
\begin{equation}\label{Pstar}
P^*(\zeta,t)=\overline{A}_0 +\sum_{j=1}^m
\frac{2\overline{A}_j\zeta}{1-\overline{\omega}_j\zeta}.
\end{equation}
On the other hand, the boundary condition on $\partial\D$ satisfied
by $\re P$ shows that
\begin{equation}\label{PPstar}
P(\zeta,t)+P^*(\zeta,t)=\frac{2q(t)}{g(\zeta,t)g^*(\zeta,t)}
\end{equation}
identically as rational functions. Therefore we find, on comparing
(\ref{P}), (\ref{Pstar}), (\ref{partialfractions}) with
(\ref{PPstar}) and using (\ref{imA0}), that $A_0$ is given by
\begin{equation}\label{A0}
A_0 =A_\infty+\sum_{j=1}^m \frac{A_j}{\omega_j}.
\end{equation}
This also shows that
\begin{equation}\label{A0mu}
A_\infty-\sum_{j=1}^m\frac{A_j}{\omega_j}=A_0-\sum_{j=1}^m\frac{2A_j}{\omega_j}=P(0)=\int
d\mu.
\end{equation}
In particular, if $m>n$, so that $A_\infty=0$, we have
\begin{equation}\label{reA0}
A_0 =\sum_{j=1}^m \frac{A_j}{\omega_j}=-\int d\mu.
\end{equation}

Now we evaluate the residues in (\ref{dynamics}):
$$
{\res}_{\zeta=\omega_k} \, \left[\zeta
P(\zeta)\frac{g'(\zeta)}{g(\zeta)}\right]
$$
$$
={\res}_{\zeta=\omega_k}  \left(A_0 \zeta+\sum_{j=1}^m
2A_j+\sum_{j=1}^m \frac{2A_j \omega_j}{\zeta-\omega_j}\right) \left(
\sum_{j=1}^m\frac{1}{\zeta-\omega_j}
-\sum_{j=1}^n\frac{1}{\zeta-\zeta_j}\right)
$$
$$
=A_0\omega_k+2A_k+\sum_{j=1, \,j\ne k}^m
\frac{2(A_k+A_j)\omega_k}{\omega_k-\omega_j}-\sum_{j=1}^n
\frac{2A_k\omega_k}{\omega_k-\zeta_j} .
$$
For the poles of $g$ the calculation of the residues is simpler
because $P(\zeta)$ is regular there: we simply have
$$
{\res}_{\zeta=\zeta_j} \, \left[\zeta
P(\zeta)\frac{g'(\zeta)}{g(\zeta)}\right] =-n_j\zeta_j P(\zeta_j),
$$
where $n_j$ is the order of the pole $\zeta_j$.

We summarize:

\begin{theorem}\label{thm:dynamics}
In the case of simple zeros of $g$ we have the rational
dynamics
\begin{equation}\label{rationalzeros}
-\frac{d}{dt}\log \omega_k=A_0+\frac{2A_k}{\omega_k}+\sum_{j=1,
\,j\ne k}^m \frac{2(A_k+A_j)}{\omega_k-\omega_j}-\sum_{j=1}^n
\frac{2A_k}{\omega_k-\zeta_j},
\end{equation}
\begin{equation}\label{rationalpoles}
-\frac{d}{dt}\log \zeta_j=A_0+\sum_{k=1}^m
\frac{2A_k}{\zeta_j-\omega_k},
\end{equation}
\begin{equation}\label{rationalb}
\frac{d\log b}{dt}=(m-n+1)A_0,
\end{equation}
where the $A_j$ are given by (\ref{C}), (\ref{Ak}), (\ref{A0}).

\end{theorem}

The last equation follows by letting $\zeta\to\infty$ in
(\ref{dotlogg}) and using (\ref{loggdot}), (\ref{loggprime}) and
(\ref{P}). The equation (\ref{rationalzeros}) was obtained, in the
polynomial case ($n=0$), already in \cite{Shraiman-Bensimon84}.

It is useful to observe that the reflected Poisson integral $P^*(\zeta,t)=\overline{P(1/\bar \zeta,t)}$ is
nothing else than the corresponding Poisson integral for the
exterior domain $\P\setminus\overline{\mathbb{D}}$:
\begin{equation}\label{poissonstar}
P^*(\zeta,t) =\frac{1}{2\pi i}
\int_{\partial(\P\setminus\overline{\mathbb{D}})}\frac{q(t)}{|g(z,t)|^2}\frac{z+\zeta}{z-\zeta}\frac{dz}{z}
\quad(\zeta\in\P\setminus\overline{\mathbb{D}}).
\end{equation}
This follows from (\ref{poisson}) and (\ref{PPstar}) by taking into
account that the jump of the Poisson integral when $\zeta$ crosses
the unit circle is exactly the right member in (\ref{PPstar}). The
orientation of the unit circle in (\ref{poissonstar}) is the
opposite of that in (\ref{poisson}). Note that $P^*>0$ in
$\P\setminus\overline{\mathbb{D}}$ with $P^*(\infty)=P(0)=\int
d\mu$, and also that the right member in (\ref{PPstar}) vanishes at
the poles $\zeta_j$ of $g$, hence that
$$
P^*(\zeta_j,t)=-P(\zeta_j,t).
$$

Next we recast parts of Theorem~\ref{thm:dynamics} more directly in
terms of the exterior Poisson integral:

\begin{theorem}\label{thm:dynamics1}
In the case of simple zeros of $g$ we have
\begin{equation}\label{rationalzerosP}
\frac{d}{dt}\log
\omega_k=P^*(\omega_k,t)-\frac{2A_k}{\omega_k}(1+\sum_{j=1}^m
\frac{1}{1-\overline{\omega}_j \omega_k}-\sum_{j=1}^n
\frac{1}{1-\overline{\zeta}_j\omega_k}),
\end{equation}
\begin{equation}\label{rationalpolesP}
\frac{d}{dt}\log \zeta_j=P^*(\zeta_j,t).
\end{equation}

For later use we recall (see (\ref{Ak})) that
\begin{equation}\label{Ak2}
\frac{2A_k}{\omega_k}=\frac{2q}{\omega_k g'(\omega_k,t)g^*(\omega_k,t)}.
\end{equation}

\end{theorem}

\begin{proof}
Set
$$
R(\zeta,t)=\frac{2q(t)}{g(\zeta,t)g^*(\zeta,t)},
$$
so that $P+P^*=R$ by (\ref{PPstar}). Returning to (\ref{dynamics})
we then have (suppressing $t$ from notation)
$$
\dot{\omega}_k =\res_{\zeta=\omega_k}[(\zeta P^*(\zeta)-\zeta
R(\zeta))(\log g(\zeta))'],
$$
where we have to compute the right member.  The residue of the first
term is immediate:
$$
\res_{\zeta=\omega_k}\zeta P^*(\zeta)(\log g(\zeta))'=\omega_k
P^*(\omega_k).
$$
For the second term, there is a slight difficulty due to the fact that both
factors in the product $R(\zeta)(\log g(\zeta))'$ have poles at
$\zeta=\omega_k$. However, we can use that
$$
(\log R)'+(\log g)'+(\log g^*)'=0
$$
by the definition of $R$, hence that
$$
\zeta R(\zeta)(\log g(\zeta))'=-\zeta{R(\zeta)}(\log
R(\zeta))'-\zeta R(\zeta)(\log g^*(\zeta))'
$$
$$
=-\zeta R'(\zeta)+\zeta^{-1}R(\zeta)(\log g(\zeta))'^*
$$
$$
=-(\zeta R(\zeta))'+R(\zeta)(1+\zeta^{-1}(\log g(\zeta))'^*)
$$
$$
=-(\zeta
R(\zeta))'+R(\zeta)(1+\sum_{j}^m\frac{1}{1-\overline{\omega}_j
\zeta}-\sum_{j}^n\frac{1}{1-\overline{\zeta}_j \zeta}).
$$
Here we used (\ref{loggprime}) in the last step.

In the last expression, the first term has no residues since it is a
pure derivative, and in the second term only $R(\zeta)$ has a pole
(a simple one) at $\omega_k$. The residue of $R(\zeta)$ is
$$
\res_{\zeta=\omega_k}R(\zeta)=
\frac{2q}{g'(\omega_k)g^*(\omega_k)}=2A_k,
$$
hence we get
$$
\res_{\zeta=\omega_k}\zeta R(\zeta)(\log g(\zeta))'=
2A_k(1+\sum_{j=1}^m
\frac{1}{1-\overline{\omega}_j\omega_k}-\sum_{j=1}^n
\frac{1}{1-\overline{\zeta}_j\omega_k}).
$$
Now (\ref{rationalzerosP}) follows. The other equation,
(\ref{rationalpolesP}), is immediate.

\end{proof}

It is easy to check that (\ref{rationalpolesP}) is valid also in the
case of multiple poles (and zeros). An immediate consequence is
that, as has previously been observed in for example \cite{Tanveer93},
\cite{Gustafsson-Prokhorov-Vasiliev04},
\cite{Gustafsson-Vasiliev06}, poles always move away from the
origin. More precisely we have

\begin{corollary}\label{cor:polesmoveout}
For each $1\leq j\leq n$ we have, when $q(t)>0$,
$$
\frac{d}{dt}\log |\zeta_j(t)| =\re {P^*(\zeta_j,t)}>0
$$
along with the estimates
\begin{equation}\label{harnack}
\frac{|\zeta_j(t)|-1}{|\zeta_j(t)|+1} \leq \frac{d\log
|\zeta_j(t)|}{d\log a_1(t)} \leq
\frac{|\zeta_j(t)|+1}{|\zeta_j(t)|-1}.
\end{equation}

\end{corollary}

\begin{proof}
It remains only to prove (\ref{harnack}), and since
$$
\frac{d\log |\zeta_j(t)|}{d\log a_1(t)}=\frac{d\log
|\zeta_j(t)|/dt}{d\log a_1(t)/dt} =\frac{\re P^*(\zeta_j,t)}{P(0)}
=\frac{\re P(\zeta_j^*,t)}{P(0)}
$$
this follows easily from the ordinary Harnack inequalities for the
positive harmonic function $\re{P(\zeta,t)}$ in $\D$, applied at the
point $\zeta=\zeta_j^*$.

\end{proof}

Because of the middle term in the right member of
(\ref{rationalzeros}), or the second term in (\ref{rationalzerosP}),
the zeros show up a more complicated behavior than the poles, and in
particular they do not always move away from the origin (see next
section). However, when $m>n$ they collectively move out, in fact even faster
than the poles. This follows by taking the
real part of (\ref{conservation}) and using (\ref{A0}),
(\ref{A0mu}), (\ref{rationalb}) and the definition (\ref{C}) of $A_\infty$:

\begin{proposition}\label{prop:conservation}
$$
\frac{d}{dt}\sum_{k=1}^m\log \omega_k(t)
-\frac{d}{dt}\sum_{j=1}^n\log \zeta_j(t)
=\begin{cases} 2(\int d\mu-A_\infty)\,\, &\text{if} \,\, m=n,\\
(m-n+2)\int d\mu\,\, &\text{if}\,\, m>n.
\end{cases}
$$
In particular, if $m>n$,
$$
\frac{d}{dt}\sum_{k=1}^m\log |\omega_k(t)|
>\frac{d}{dt}\sum_{j=1}^n\log |\zeta_j(t)|,
$$
$$
\frac{d}{dt}\sum_{k=1}^m\arg\omega_k(t)
=\frac{d}{dt}\sum_{j=1}^n\arg\zeta_j(t).
$$
\end{proposition}

The case $m=n$ is indeed a little exceptional in that $f(\zeta,t)$
in this case has a simple pole at infinity which may be lost
at one moment of time (see discussion after (\ref{structureg})), causing $m$ to
temporarily drop below $n$. What happens in (\ref{structureg}) then is that (at least)
one root $\omega_k(t)$ rushes to infinity while $b(t)$ goes to zero.
Then also $A_\infty$ goes to infinity. However, after the event the root comes back again.

\begin{example}\label{ex:drop}
The following example of an  off-center injection/suction of a disk
is similar to examples which have been given by Y.~Hohlov, S.~Howison, S.~Richardson and others,
see \cite{Hohlov-Howison93}, \cite{Richardson97}.
We consider a Hele-Shaw evolution with $q(t)>0$ such that at one particular
instant, which we choose to be $t=2$, the fluid domain is $\Omega(2)=\D(1,2)$.
Then for a general $t$ the domain $\Omega(t)$ will satisfy the quadrature identity
\begin{equation}\label{twoptqi}
\int_{\Omega(t)}hdm = 2\pi (Q(t)-Q(2))h(0)+ 4\pi h(1)
\end{equation}
for functions $h$ harmonic and integrable in the domain.
The mapping function at $t=2$ is
$$
f(\zeta,2)=\frac{3\zeta}{2-\zeta},
$$
while it for $t\ne 2$ can be taken to be of the form
$$
f(\zeta,t)=b(t)\frac{\zeta(\zeta-a(t))}{\zeta -t}
$$
for suitable real-valued functions $a(t)$ and $b(t)$. Here we have left the
exact choice of $q(t)$ open in order to allow the time coordinate simply
to be the location of the pole on the positive real axis, which by
Corollary~\ref{cor:polesmoveout} is consistent with the assumption
$q(t)>0$. The range of $t$ will then be $1<t<\infty$.

The derivative of $f$ is
$$
g(\zeta,t)=b(t)\frac{\zeta^2-2t\zeta +ta(t)}{(\zeta-t)^2},
$$
which has two zeros $\omega_1(t)$, $\omega_2(t)$ satisfying
\begin{equation}\label{twozeros}
\begin{cases}
\omega_1 \omega_2 =ta(t),\\
\omega_1 +\omega_2 =2t.
\end{cases}
\end{equation}
Using the fixed data, contained in the last term in (\ref{twoptqi}), one gets the equations
$$
f(t^{-1},t)=\frac{b(1-at)}{t(1-t^2)}=1,
$$
$$
\res_{\zeta=1/t} f^*(\zeta,t)g(\zeta,t)d\zeta =\frac{b^2(1-2t^2+at^3)(a-t)}{t(1-t^2)^2}=4,
$$
which can be solved for $a$, $b$, giving
$$
a(t)=\frac{1}{t}+\frac{(t^2-1)(t^2+1+\sqrt{(t^2-1)^2+16})}{2t(t^2-4)},
$$
$$
b(t)=\frac{t}{2}(t^2+1-\sqrt{(t^2-1)^2+16}).
$$
Now the behavior of the roots can be read off from (\ref{twozeros}):
for $1<t<2$ we have $\omega_1\omega_2=at<0$, tending to $-\infty$ as
$t\to 2$, while $\omega_1+\omega_2=2t$ remains finite. Also, $b(t)\to 0$ as $t\to 2$.
Thus the two roots go to infinity along the real axis (in different directions).

For $2<t<\infty$, $\omega_1\omega_2=at>0$, which means that the roots come back along the imaginary
axis as $t$ increases from $2$. Finally, as $t\to \infty$, they turn back to infinity
the same way.

\end{example}


\section{Global locally univalent solutions are univalent}
\label{sec:locallyunivalent}

In this section we show that any global locally univalent solution of
the Polubarinova-Galin equation actually has to be univalent all
the time. It is intuitively clear that it has to be so. In fact, if
univalence breaks down while local univalence still holds, then this
means that two different parts of $f(\mathbb{D},t)$ start to
overlap, hence a hole in the part of $\C$ covered by
$f(\mathbb{D},t)$ is created. With the solution being global this
hole eventually has to be filled in, but it is easy to realize that
when this occurs a branch point for the covering map from the
multiply sheeted $f(\mathbb{D},t)$ to $\mathbb{C}$ is created. In
other words, a zero of $f'$ moves into $\D$ and local univalence is
lost.

A rigorous proof of the above statement can be based on a result
from L\"owner theory, namely Theorem 6.2 in \cite{Pommerenke75}. To prepare
for that, recall that in terms of the Taylor expansion
(\ref{powerseries})
of $f(\zeta,t)$ the zeroth order moment is given by
$$
M_0(t)=\sum_{j=1}^{\infty}j|a_{j}(t)|^{2}.
$$
Let
\begin{equation}\label{N0}
N_0(t)=\sum_{j=2}^{\infty}j|a_{j}(t)|^{2}
\end{equation}
denote what remains after the first term. Then
$$
N_0(t)=M_0(t)-a_1(t)^2=M_0(0)+2Q(t)-a_1(t)^2
$$
(recall (\ref{Q})).

In \cite{Kuznetsova01} O.~Kuznetsova showed that $N_0(t)$ is a decreasing
function of $t$, hence that $a_{1}(t)^2-2Q(t)$ is increasing and
that
\begin{equation}\label{diffu}
0\leq N_0(t)\leq N_0(0),
\end{equation}
\begin{equation}\label{diffusecond}
a_{1}(t)^2\geq a_{1}(0)^2+2Q(t).
\end{equation}
It also follows that
\begin{equation}\label{diffu1}
\frac{a_1(t)}{\sqrt{M_0(0)+2Q(t)}}=\frac{a_1(t)}{\sqrt{M_0(t)}}\nearrow
1 \quad {\rm as\,\,} t\to\infty,
\end{equation}
and that
\begin{equation}\label{akN}
|a_k(t)|\leq \sqrt{\frac{N_0(0)}{k}} \quad {\rm for\,\,} k\geq 2.
\end{equation}

\begin{theorem}\label{lemma2}
Let $f(\zeta,t)\in\mathcal{O}_{\rm locu}(\overline{\mathbb{D}})$ be
a global ($0\leq t<\infty$) locally  univalent solution of the
Polubarinova-Galin equation (\ref{pg1}) driven by injection with
rate $q(t)>0$ such that $Q(t)\to\infty$ as $t\to \infty$. Then $
f(\zeta,t)\in \mathcal{O}_{\rm univ}(\overline{\mathbb{D}}) $ for
all $0\leq t<\infty$.
\end{theorem}

\begin{proof}
The choice of time parameter is irrelevant as long as $Q(t)\to
\infty$ as $t\to\infty$. Therefore, by (\ref{diffu1}), we may choose
$t$ so that
$$
a_1(t)=e^t.
$$
Then $f(\zeta,t)$ solves the L\"owner-Kufarev equation (\ref{lk})
with $P(\zeta,t)$ being a Carath\'{e}odory function, i.e.,
satisfying
$$
{\rm Re\,} P(\zeta,t)>0 \quad {\rm in\,\,} \mathbb{D},\quad
P(0,t)=1.
$$
This is the setting in \cite{Pommerenke75}.

Now, for (say) $|\zeta|\leq \frac{1}{2}$ and $t\geq 0$ we have,
using (\ref{diffu}),
\begin{align}
|f(\zeta,t)|&\leq  e^{t}+\sum_{j=2}^{\infty}|a_{j}(t)|2^{-j}\notag\\
       &\leq e^{t}+\left(\sum_{j=2}^{\infty}j|a_{j}({t})|^2\right)^{1/2}
          \left(\sum_{j=2}^{\infty}\frac{1}{j}2^{-2j}\right)^{1/2}\notag\\
          &\leq e^t +\sqrt{N_0(t)}\leq  e^t +\sqrt{N_0(0)}.\notag
\end{align}
Thus $|f(\zeta,t)|\leq Ce^t$ ($0\leq t<\infty$) in a neighborhood of
the origin, which means that the assumptions in Theorem 6.2 of
\cite{Pommerenke75} are satisfied. The conclusion is that $f(\zeta,t)$ is
univalent for all $t\geq 0$. In our case it follows that $f(\zeta,t)$ is
actually univalent in the closed unit disk, because by assumption it
is locally univalent there, and if for two different points
$\zeta_1, \zeta_2\in\partial \mathbb{D}$ we had
$f(\zeta_1,t)=f(\zeta_2,t)$ for some $t$ then univalence in the open
disk $\mathbb{D}$ would be lost in the next instance (the boundary
is really propagating with positive speed under our assumptions).
\end{proof}

\begin{remark}
Under the slightly weaker assumption that
$f(\zeta,t)\in\mathcal{O}_{\rm norm}(\overline{\mathbb{D}})\cap\mathcal{O}_{\rm locu}({\mathbb{D}})$
(i.e., local univalence only in the open disk)
the theorem holds with the conclusion correspondingly changed to
$f(\zeta,t)\in\mathcal{O}_{\rm univ}({\mathbb{D}})$.
\end{remark}


\section{Asymptotic questions in the polynomial case}\label{sec:asymptotic}

In this section we study asymptotic behaviour of the roots of $g=f'$
for locally univalent polynomial solutions of (\ref{pg1}) with
$q(t)>0$ and $Q(t)\to \infty$ as $t\to \infty$. Thus $n=0$  and we
write
\begin{equation}\label{polynomial}
f(\zeta,t)=\sum_{j=1}^{m+1}a_{j}(t)\zeta^{j}
\end{equation}
We assume $a_{m+1}(0)\ne 0$, by which $a_{m+1}(t)\ne 0$ for all $t$,
since $M_m=\overline{a_{m+1}(t)}a_1(t)^{m+1}$ (see
(\ref{lastmoment}) below) is preserved.


\subsection{Collision of roots}\label{sec:basicquestions}

In the rational case we have seen (Example~\ref{ex:drop}) that the zeros of $g$ may occasionally
(when $m=n$) move to infinity and come back again, and we have seen that the poles always
move out. In the polynomial case ($n=0$) the zeros cannot move to infinity in finite time,
as can be seen from Proposition~\ref{prop:conservation} for example. Below we give some examples showing that
zeros can still move towards the origin and that they can collide, even though they in principle
repel each other (as can be seen from (\ref{rationalzeros}) for example).

\begin{example}\label{ex:collision}
When $m=2$, $n=0$ and the mapping function $f$ has
real coefficients the solution of the Polubarinova-Galin equation
can be made fully explicit (see Example~\ref{ex:huntingford}). The roots $\omega_1$ and
$\omega_2$ will either be a complex conjugate pair, or both will be
real, with occasional collisions allowed. We shall first make some
immediate conclusions from the dynamical equations
(\ref{rationalzeros}) in the case that the roots $\omega_1$,
$\omega_2$ are real.

Assume for example that both roots are positive, say
$$
1<\omega_1 <\omega_2.
$$
By (\ref{Ak}) the coefficients $A_k$ are then given by
$$
\frac{|b|^2 }{q}A_1 =\frac{\omega_1^2}{(\omega_1
-\omega_2)(1-\omega_1^2)(1-\omega_1\omega_2)},
$$
$$
\frac{|b|^2 }{q}A_2 =\frac{\omega_2^2}{(\omega_2
-\omega_1)(1-\omega_2^2)(1-\omega_1\omega_2)}.
$$
Now (\ref{rationalzeros}) together with (\ref{A0}) gives
$$
\frac{|b|^2}{q}\frac{d}{dt}\log \omega_1 =\alpha (\omega_1,\omega_2)
+\beta (\omega_1,\omega_2)
\cdot\frac{\omega_1\omega_2-3}{\omega_2-\omega_1},
$$
$$
\frac{|b|^2}{q}\frac{d}{dt}\log \omega_2 =\alpha (\omega_1,\omega_2)
-\beta(\omega_1,\omega_2)
\cdot\frac{\omega_1\omega_2-3}{\omega_2-\omega_1},
$$
where
$$
\alpha (\omega_1,\omega_2)=\frac{2(1+\omega_1\omega_2)}{(\omega_1^2-1)(\omega_2^2-1)(\omega_1\omega_2-1)},
$$
$$
\beta(\omega_1,\omega_2)=\frac{(\omega_1+\omega_2)}{(\omega_1^2-1)(\omega_2^2-1)(\omega_1\omega_2-1)}.
$$
Note that $\alpha (\omega_1,\omega_2)$ and $\beta(\omega_1,\omega_2)$ are
positive and symmetric. The remaining factor
$\frac{\omega_1\omega_2-3}{\omega_2-\omega_1}$ has a singularity
when $\omega_1=\omega_2$ and changes sign when $\omega_1\omega_2=3$.
Hence we can draw the following conclusions.

\begin{itemize}

\item The product $\omega_1\omega_2$ always increases in time (also clear from
Proposition~\ref{prop:conservation}).

\item The ratio $\omega_2/\omega_1$ decreases if
$\omega_1\omega_2<{3}$, increases if $\omega_1\omega_2>{3}$

\item  If $\omega_1\omega_2<{3}$ and $\omega_2-\omega_1$ is sufficiently small, then $\dot\omega_1<0$
and $\dot\omega_2>0$. Hence no collisions occur for
$\omega_1\omega_2<3$.

\item  If $\omega_1\omega_2={3}$ then $\dot\omega_1=\dot\omega_2>0$.

\item  If $\omega_1\omega_2>{3}$ and $\omega_2-\omega_1$ is sufficiently small, then $\dot\omega_1>0$
and $\dot\omega_2<0$. The condition $\omega_1\omega_2>{3}$ is
preserved in time and eventually leads to a collision (even if
$\omega_2-\omega_1$ is not small to start with). After the collision
the roots cease to be real and continue as a complex conjugate pair.

\end{itemize}

\end{example}

In the case of complex roots (but still with $f$ having real
coefficients) the scenario may be quite interesting, as  the
following example due to C.~Huntingford \cite{Huntingford95} shows.

\begin{example}\label{ex:huntingford}
When expressed in terms of the moments $M_1$, $M_2$ and with the
time parameter chosen so that $a_1(t)=e^t$ ($q(t)>0$ defined
accordingly), the solution of the Polubarinova-Galin equation
(\ref{pg1}) in the case $m=2$, $n=0$ with $f$ having  real
coefficients is explicitly
$$
f(\zeta,t) = e^t \zeta +\frac{M_1\zeta^2}{e^{2t}+3e^{-2t}M_2} +
e^{-3t}M_2 \zeta^3.
$$

By taking
$$
\begin{cases}
M_1=\frac{32}{25}, \\
M_2=\frac{1}{5} \end{cases}
$$
one gets an interesting example which has been
investigated by C.~Huntingford \cite{Huntingford95}.
We let the solution start at the first possible instant of time,
namely $t=t_0=\frac{1}{4}\log\frac{3}{5}<0$ (for smaller values of
$t$, $f$ is not locally univalent). At that moment the image domain
$\Omega(t)=f(\D,t)$ has two cusps on the boundary (the two zeros of
$g$ form a complex conjugate pair on the unit circle). As $t$
increases these cusps resolve, the zeros of $g$ move out from the
unit circle and collide on the real axis at some moment $t=t_1<0$.
After that one of the zeros, say $\omega_1(t)$, moves back to the
unit circle and reaches it again at time $t=0$. Thus a new cusp
(which will be a $5/2$-power cusp) develops on the boundary. However
also this cusp resolves, the root $\omega_1(t)$ moves away from the
unit circle along the real axis, captures and collides with the
other root $\omega_2(t)$, after which they leave the real axis and
finally move towards infinity in the asymptotic directions of the
positive and negative imaginary axes. All of this can be seen by
examining \cite{Huntingford95} carefully. The final asymptotics
follows from Theorem~\ref{thmgg} below.

The solution is global ($t_0<t<\infty$) and is all the time
univalent (otherwise it would have broken down in finite time, by
Theorem~\ref{lemma2}), or the remark following it.

\end{example}

The above example shows that it is possible for a single zero to
reach the unit circle in the injection case.  The next example will
show that this can not occur if the other zeros are sufficiently
separated and far away.

\begin{example}\label{ex:onecloseroot}
Let $n=0$, $m\geq 2$ with $f$ allowed to have complex coefficients.
We assume that one of the roots is close to the unit circle and the
others far away and well separated from each other, say
$$
1<\omega_{1}<1+\varepsilon, \quad \varepsilon>0 \,\,{\rm small},
$$
$$
|\omega_{k}-\omega_{j}|\geq M\quad (1\leq k, j\leq m,\,\, k\neq
j,\,\, M>1\,\,\text{large}).
$$
We shall use equation~(\ref{rationalzeros}) to investigate short
time root behavior.

First observe by (\ref{Ak}) that
$$
\frac{A_k}{\omega_k}=\frac{q}{|b|^2(1-|\omega_k|^2)}\cdot\frac{1}
{\prod_{j\neq
k}(\omega_{k}-\omega_{j})\overline{(\omega_{k}^{*}-\omega_{j})}},
$$
hence that, among these quantities, $\frac{A_1}{\omega_1}$ is the
dominating one under the present assumptions. Therefore, by
(\ref{rationalzeros}) and (\ref{A0}) (with $A_\infty=0$) we have
$$
-\frac{d}{dt}\log|\omega_{1}(t)| =\re{A_0}+\re{\frac{2A_1}{\omega_1}}
+\re\sum_{j=2}^m\frac{2(A_1+A_j)}{\omega_1-\omega_j}
$$
$$
\approx
\re\frac{3A_1}{\omega_1}+\re[\frac{2A_1}{\omega_1}\sum_{j=2}^m\frac{\omega_1}{\omega_1-\omega_j}]
$$
$$
\approx
\re\frac{3A_1}{\omega_1}
=\frac{3q}{|b|^2(1-|\omega_1|^2)}\cdot\re\frac{1}
{\prod_{j=2}^m(\omega_{1}-\omega_{j})\overline{(\omega_{1}^{*}-\omega_{j})}}
$$
$$
\approx\frac{3q}{|b|^2(1-|\omega_1|^2)\prod_{j=2}^m|\omega_{1}-\omega_{j}|^2}<0,
$$
which shows that the minimum root $\omega_1(t)$ moves away from the
origin to prevent the strong solution from blowing up in short time.

Similarly, in the case of suction, under the same assumptions the
minimum root $\omega_{1}$ moves towards the origin initially and will make the solution blow up.
\end{example}


\subsection{Long time behaviour of roots and coefficients}\label{sec:sufficient}

We now assume that the locally univalent polynomial solution
(\ref{polynomial}) is global in time. By Theorem~\ref{lemma2},
$f(\zeta,t)$ is then actually univalent all the time. We shall
describe the asymptotic behavior of roots of $g(\zeta,t)$. It will
be shown (Theorem~\ref{thmgg}) that these roots all move away from
the origin as time $t$ is large, even though some of them may move
towards the origin initially. One naturally expects that what makes the
initial root behavior of $g(\zeta,t)$ unpredictable is that the
distribution of zeros of $g(\zeta,0)$ is not always uniform.
Therefore, it seems reasonable that, by controlling the distribution
of these initial roots, we can guarantee that all roots always move
out. In fact, this will be demonstrated in Theorem~\ref{tem}.

Recall first Richardson's formula \cite{Richardson72} for the
harmonic moments $M_{k}$ (see (\ref{defmoments})):
\begin{equation}\label{moments}
M_{k}=\sum_{(i_{1},\cdots,i_{k+1})}i_{1}a_{i_{1}}\cdots
a_{i_{k+1}}\overline{a_{i_{1} +\cdots+i_{k+1}}}.
\end{equation}
Here $k=0,1,2,\dots$ and the summation runs over all $(k+1)$-tuples
$({i_{1},\cdots,i_{k+1}})$ of integers with $i_j\geq 1$, with the
convention that $a_j=0$ for $j>m+1$ in the present polynomial case
(\ref{polynomial}). See \cite{Richardson72} and \cite{Kuznetsova01}
for more details. Note that the final nonzero moment has a quite
simple expression:
\begin{equation}\label{lastmoment}
M_m=a_1^{m+1}\overline{a_{m+1}}.
\end{equation}
Lemma~\ref{decay} below will show that asymptotically, as
$t\to\infty$, the formula (\ref{lastmoment}) will almost be true for
all moments.

Besides the zeros $\omega_{k}(t)$ of $g(\zeta,t)$, which are what we
want to study, we introduce two other sets of zeros, for comparison:

\begin{itemize}

\item Let
\begin{equation}\label{omegaktilde}
\tilde{\omega}_{k}(t)={\omega_{k}(t)}{a_{1}(t)^{-\frac{m+2}{m}}}
\end{equation}
be corresponding {\bf rescaled zeros}, namely the zeros of
$\tilde{g}(\zeta,t)=g(a_1^{\frac{m+2}{m}}\zeta,t)$.

\item Let
$$
\hat{\omega}_k=\frac{1}{\sqrt[m]{-{(m+1)\overline{M}_{m}}}},
$$
be the zeros of the polynomial $a_1+(m+1)a_1^{m+2}a_{m+1}\zeta^m$
consisting of the first and last term in $\tilde{g}(\zeta,t)$.

\end{itemize}

Note that the $\hat{\omega}_k$ do not depend on time. On the other hand,
the rescaled zeros $\tilde{\omega}_{k}(t)$ are normalized in such a
way that they asymptotically stabilize, and approach the
$\hat{\omega}_k$ (Theorem~\ref{thmgg}).
Note also that the product $\tilde{\omega}_1(t)\cdots
\tilde{\omega}_m(t)=-\frac{1}{(m+1)\overline{M}_m}$ is a conserved quantity.

The following lemma is a slightly sharpened version of estimates first
obtained in \cite{Lin09a}.

\begin{lemma}\label{decay}
Given $m\geq 1$ there exist universal polynomials
$F_j(x_1, \dots, x_{m+2})$, $1\leq j\leq m+1$, in $m+2$ real variables such that
\begin{equation}\label{moment0}
\left|a_{j}(t)a_{1}(t)^{j}-\overline{M_{j-1}}\right| \leq
\frac{1}{a_{1}(t)^4} F_{j}(\frac{1}{a_{1}(0)},\sqrt{N_0(0)},
|M_{1}|,\cdots,|M_{m}|)
\end{equation}
($t>0$, $1\leq j\leq m+1$) whenever $f(\zeta,t)$ given by
(\ref{polynomial}) is a locally univalent solution of (\ref{pg1}).
The polynomials $F_{j}$ are
increasing functions in all their variables.

In particular we have estimates of the kind
\begin{equation}\label{estimateN0}
|N_0(t)|\leq \frac{C_1}{a_{1}(t)^4},
\end{equation}
\begin{equation}\label{aj}
|a_j(t)|\leq \frac{C_j}{a_1(t)^j}, \quad 2\leq j \leq m+1,
\end{equation}
for suitable constants $C_j=C_j(\frac{1}{a_{1}(0)},\sqrt{N_0(0)},
|M_{1}|,|M_{2}|,\cdots,|M_{m}|)$.
\end{lemma}

\begin{proof}
Using (\ref{moments}) we first prove by induction on decreasing $j$
that there exist polynomials $G_j(x_1,\dots,x_{m+2})$, $1\leq j\leq m+1$, in $m+2$
variables, increasing in all variables, such that
\begin{equation}\label{moment1}
\left|a_{j}(t)a_{1}(t)^{j}-\overline{M_{j-1}}\right|\leq
\frac{1}{a_{1}(t)^2} G_{j}(\frac{1}{a_{1}(0)}, \sqrt{N_0(0)}, |M_{1}|,\cdots,
|M_{m}|).
\end{equation}

We start by setting $G_{m+1}=0$, which makes (\ref{moment1}) hold
for $j=m+1$. Now we show (\ref{moment1}) for $j=k$ assuming it holds
for all $j\geq k+1$. The arguments will work for any $1\leq k \leq
m$. Using (\ref{moments}) and suppressing dependence on $t$ from
notation, we have
$$
{M_{k-1}}-\bar{a}_k a_1^k=\sum_{(i_{1},\cdots,i_{k})\ne
(1,\dots,1)}i_{1}a_{i_{1}}\cdots a_{i_{k}}\cdot\overline{a_{i_{1}
+\cdots+i_{k}}}
$$
$$
=\sum_{(i_{1},\cdots,i_{k})\ne
(1,\dots,1)}\frac{i_{1}a_{i_{1}}\cdots a_{i_{k}}}{a_1^{i_1+\dots
+i_k}}\cdot \overline{a_{i_{1} +\cdots+i_{k}}} a_1^{i_1+\dots +i_k}
$$
$$
=\frac{1}{a_1^2}\sum_{(i_{1},\cdots,i_{k})\ne
(1,\dots,1)}\frac{i_{1}a_{i_{1}}\cdots a_{i_{k}}}{a_1^{i_1+\dots
+i_k-2}} \cdot\overline{a_{i_{1} +\cdots+i_{k}} a_1^{i_1+\dots
+i_k}}.
$$

We shall estimate the terms in the above sum. Since
${(i_{1},\cdots,i_{k})\ne (1,\dots,1)}$, at least one of the $i_j$ is
$\geq 2$, hence it follows from (\ref{akN}), (\ref{diffusecond}) that the
first factors can be estimated as
\begin{equation}\label{estimate}
\left|\frac{i_{1}a_{i_{1}}\cdots a_{i_{k}}}{a_1^{i_1+\dots +i_k-2}}\right| \leq
(m+1)\sum_{0\leq \alpha,\beta\leq
m+1}\frac{(\sqrt{N_0(0)})^\alpha}{a_1(0)^\beta}.
\end{equation}
Moreover $i_1+\dots + i_k \geq k+1$, so by
the induction hypothesis we have
$$
|{a_{i_{1} +\cdots+i_{k}} a_1^{i_1+\dots
+i_k}}-\overline{M_{i_1+\dots +i_k-1}}|
$$
$$
\leq G_{i_1+\dots +i_k} (\frac{1}{a_{1}(0)}, \sqrt{N_0(0)},
|M_{1}|,\cdots, |M_{m}|)\frac{1}{a_{1}(t)^2},
$$
in particular
$$
|{a_{i_{1} +\cdots+i_{k}} a_1^{i_1+\dots +i_k}}|
$$
$$
\leq G_{i_1+\dots +i_k} (\frac{1}{a_{1}(0)}, \sqrt{N_0(0)},
|M_{1}|,\cdots, |M_{m}|)\frac{1}{a_{1}(0)^2} +\sum_{j=1}^m |M_j|.
$$
Therefore we can estimate also the second factors in the above
expression for $a_ka_1^k-\overline{M_{k-1}}$. From this we easily
deduce (\ref{moment1}) for $j=k$ knowing that it is true for $j\geq
k+1$. This completes the induction step and hence proves
(\ref{moment1}) for $1\leq j \leq m$.

Now, already (\ref{moment1}) shows that the estimate (\ref{aj})
holds. Since this estimate improves (\ref{akN}) by a factor
at least $a_1(t)^2$ in the denominator, we can `bootstrap' the
previous argument: using the new estimate in (\ref{estimate}) makes
the induction process work with the factor $1/a_1(t)^2$ in
(\ref{moment1}) replaced by $1/a_1(t)^4$. Thus (\ref{moment0})
follows, and since (\ref{estimateN0}) is just the special case $j=1$
the lemma is proved.

\end{proof}

\begin{theorem}\label{thmgg}
With notations as above, assume that $f(\zeta,t)$ is a
global polynomial solution of (\ref{pg1}). Then
\[\left|\tilde{\omega}_{k}(t)-\hat{\omega}_{k}\right|\rightarrow 0
\quad \mbox{as $t$}\rightarrow\infty\] if the roots
$\tilde{\omega}_{k}(t)$ and $\hat{\omega}_{k}$ are ordered
appropriately. Furthermore, all roots eventually move away from the
origin as time is large enough.
\end{theorem}

\begin{proof}
The monic polynomial vanishing at the rescaled roots
(\ref{omegaktilde}) is
\begin{align*}
\prod_{j=1}^{m}(\zeta-\tilde{\omega}_{j})&
=\sum_{j=1}^{m+1}\left({a_{1}^{-\frac{m+2}{m}}}\right)^{m+1-j}
\frac{ja_{j}}{(m+1)a_{m+1}}\,\zeta^{j-1}\notag\\
&=\sum_{j=1}^{m+1}\frac{j}{m+1}\,\frac{a_{j}a_{1}^{j}}
{a_{m+1}a_{1}^{m+1}}a_{1}^{-\frac{2}{m}(m+1-j)}\,\zeta^{j-1}\notag\\
&=\frac{1}{(m+1)\overline{M_{m}}}+\sum_{j=2}^{m}\frac{j}{m+1}\,\frac{a_{j}a_{1}^{j}}
{\overline{M_{m}}}a_{1}^{-\frac{2}{m}(m+1-j)}\,\zeta^{j-1}+\zeta^{m}.
\end{align*}
Due to Lemma~\ref{decay} and because the exponents
${-\frac{2}{m}(m+1-j)}$ are strictly negative the coefficients of
the middle terms tend to zero:
$$
\left|\frac{j}{m+1}\,\frac{a_{j}a_{1}^{j}}
{\overline{M_{m}}}a_{1}^{-\frac{2}{m}(m+1-j)}\right|
\leq\left|\frac{j}{m+1}\,\frac{F_j}
{a_{1}^{4}{M_{m}}}a_{1}^{-\frac{2}{m}(m+1-j)}\right|
$$
\begin{equation}\label{rescaling}
+\left|\frac{j}{m+1}\,\frac{M_{j-1}}
{{M_{m}}}a_{1}^{-\frac{2}{m}(m+1-j)}\right|\to 0 \quad {\rm
as\,\,} t\to \infty \quad (2\leq j\leq m).
\end{equation}
From this it follows that
$|\tilde{\omega}_{k}(t)-\hat{\omega}_{k}|\rightarrow 0$ as
$t\rightarrow\infty$.

It also follows that the roots move away from the origin as time is
large, because the speed of the roots only depends on the position of
the roots, and it is clear that for a symmetric configuration of
roots, like $\hat{\omega}_k$, the speed points radially away from the
origin. A slightly more precise argument can be based on
(\ref{rationalzerosP}), (\ref{Ak2}), by which
$$
\frac{d}{dt}\log |\omega_k|=\re P^*(\omega_k)-\re\frac{2q}{\omega_k
g'(\omega_k)g^*(\omega_k)}(1+\sum_{j=1}^m
\frac{1}{1-\overline{\omega}_j\omega_k}).
$$
Here the first term is always positive, while the subtracted term is
negative whenever the roots are sufficiently far away and close to
symmetrical configuration. In
fact, in this case $\omega_{k}g'(\omega_{k})\approx -ma_{1}$,
$g^{*}({\omega_{k}})\approx a_{1}$ and
$1+\sum_{j=1}^{m}\frac{1}{1-\overline{\omega_{j}}\omega_{k}}\approx
1$. (This argument will be made more precise in the proof of
Theorem~\ref{tem} below.) Hence $\frac{d}{dt}\log |\omega_k|>0$
for large $t$, as claimed.

\end{proof}


\subsection{Assumptions only on initial data}

In this subsection we do not assume {\it a priori} that $f(\zeta,t)$ is
global, we only make assumption on the initial data, and the
solution being global will be part of the conclusion.

\begin{theorem}\label{tem}
In terms of any initial function
$f(\zeta,0)=\sum_{j=1}^{m+1}a_{j}(0)\zeta^{j}$, let $M>0$ be a
common upper bound for $|M_{1}|, \dots, |M_{m}|$ and $\sqrt{N_0(0)}$.
Then, for any $\varepsilon>0$, there exists a number $B=B(\varepsilon, M, m)>0$ such
that whenever $a_{1}(0)>B$ the polynomial solution $f(\zeta,t)$
starting with $f(\zeta,0)$ is global in time and the assertions (i)-(iii)
below hold for $1\leq k\leq m$.

\begin{itemize}

\item[(i)] For a suitable ordering of the $\tilde{\omega}_{k}(t)$ and $\hat{\omega}_{k}$
we have
$$
\left|\tilde{\omega}_{k}(t)-\hat{\omega}_{k}\right|<\varepsilon,\quad 0\leq
t<\infty.
$$
In particular the roots $\omega_{k}$ never collide (with $\varepsilon >0$
sufficiently small).

\item[(ii)]
$$
\left|\tilde{\omega}_{k}(t)-\hat{\omega}_{k}\right|\rightarrow 0
\quad\mbox{as $t\rightarrow\infty$}.
$$

\item[(iii)] The roots $\omega_{k}(t)$ move away from the origin (for all $t\geq 0$).

\end{itemize}

\end{theorem}

\begin{proof}
First, choose  $\rho(M_{m}, m)>0$ so that the disks
$\D(\hat{\omega}_{j},\rho)$ ($1\leq j\leq m$) are disjoint and such
that, for any point $z\in\D(\hat{\omega}_{j},\rho)$,
\begin{equation}\label{delta}
\left|z^{m}+\frac{1}{(m+1)\overline{M_{m}}}\right|
<\frac{1}{4(m+1)|M_{m}|}\tan\frac{\pi}{20}.
\end{equation}
For example, any $0<\rho<
\frac{1}{4}|M_m|^{-\frac{1}{m}}\tan\frac{\pi}{20}$ will do.

Next, as in the proof of Theorem~\ref{thmgg} we have that for any $0<\varepsilon<\rho$
there exists $\delta>0$ such that $|\tilde{\omega}_k-\hat{\omega}_k|<\varepsilon$ whenever
\begin{equation}\label{middleterm}
\left|\frac{j}{m+1}\frac{a_{j}a_{1}^{j}}
{\overline{M_{m}}}a_{1}^{-\frac{2}{m}(m+1-j)}\right|<\delta,
\quad 2\leq j\leq m.
\end{equation}
Finally, it follows from Lemma~\ref{decay} that (\ref{middleterm}) indeed holds if just
$a_1(0)$ is large enough.
We conclude that for any $0<\varepsilon<\rho$, there exists $B_{1}(\varepsilon, M, m)>0$
such that if $a_{1}(0)>B_{1}$, then
$|\tilde{\omega}_{k}(t)-\hat{\omega}_{k}|<\varepsilon$
and (\ref{delta}) holds.  Therefore $(i)$ holds. By
Theorem~\ref{thmgg}, also $(ii)$ holds.

From (\ref{rationalzerosP}), (\ref{Ak2}) we see that in order to prove that all
roots move away from the origin, it is sufficient to prove that
\begin{equation}\label{ggo}
\re\left[\frac{2}{\omega_{k}g'(\omega_{k},t)g^{*}(\omega_{k},t)}
\left(1+\sum_{j=1}^{m}\frac{1}{1-\overline{\omega_{j}}\omega_{k}}\right)\right]<0.
\end{equation}
We first estimate the terms of $\omega_{k}g'(\omega_{k},t)
=\sum_{j=1}^{m}j(j+1)a_{j+1}\omega_{k}^{j}$. For $1\leq j\leq m-1$
we have
\begin{align*}
a_{j+1}\omega_{k}^{j}&=a_{j+1}\left(\frac{\omega_{k}}
{a_{1}^{\frac{m+2}{m}}}\right)^{j}a_{1}^{\frac{m+2}{m}j}\notag\\
&=\left(\frac{\omega_{k}}{a_{1}^{\frac{m+2}{m}}}\right)^{j}
\cdot a_{1}^{j+1}a_{j+1}\cdot a_{1}^{\frac{m+2}{m}j-(j+1)}\notag\\
&=\tilde{\omega}_{k}^{j}\cdot a_{1}^{j+1}a_{j+1}\cdot
a_{1}^{2(\frac{j}{m}-1)}a_{1}.
\end{align*}
Hence, for $1\leq j\leq m-1$,
\begin{equation}\label{com}
\left|a_{j+1}\omega_{k}^{j}\right|
\leq\left[1+\left(\varepsilon+|\hat{\omega}_{k}|\right)^{m-1}\right]
\left(M+ \frac{F_j}{a_1^4}\right)a_{1}^{-2/m}a_{1}
\end{equation}
with $F_{j}=F_{j}(\frac{1}{a_{1}(0)}, M,\cdots,
M)$ as in Lemma~\ref{decay}.
From (\ref{com})
and (\ref{diffusecond}) we conclude that there exists $B_{2}>B_1$
such that if $a_{1}(0)>B_2$, then
$|j(j+1)a_{j+1} \omega_{k}^{j}|\leq\frac{1}{4}a_{1}\tan\frac{\pi}{20}$.
For $j=m$, $a_{j+1}\omega_k^j=a_{m+1}\omega_k^m=\tilde{\omega}_k^m \overline{M}_j a_1$
by (\ref{omegaktilde}).

When $a_{1}(0)>B_{2}$ we have
\begin{align*}
\omega_{k}g'(\omega_{k},t)+ma_{1}&=\left(\tilde{\omega}_{k}^{m}
+\frac{1}{\overline{M}_{m}(m+1)}\right)\overline{M_{m}}m(m+1)a_{1}\notag\\
&+\sum_{j=1}^{m-1}j(j+1)a_{j+1}\omega_{k}^{j}
\end{align*}
due to (\ref{delta}), and since
$|j(j+1)a_{j+1}\omega_{k}^{j}|\leq\frac{1}{4}a_{1}\tan\frac{\pi}{20}$ we conclude from this that $|\arg\omega_{k}(t)g'(\omega_{k},t)-\pi|<\frac{\pi}{20}$. Finally,  we can find
$B=B(\delta, M, m)>B_{2}$ such that if $a_{1}(0)>B$, then
$|\arg \omega_{k}(t)g'(\omega_{k},t)-\pi|$, $|\arg
g^{*}({\omega_{k}},t)|$ and
$|\arg(1+\sum_{j=1}^{m}\frac{1}{\overline{\omega_{j}}\omega_{k}})|$
are all $<\frac{\pi}{20}$. Therefore (\ref{ggo}), and hence $(iii)$ in the theorem, holds.

\end{proof}



\section{Asymptotics for rational solutions}
\label{sec:multiplecut}

In this section we study the asymptotics of poles and Taylor coefficients in the
rational case. We assume that $f\in\mathcal{O}_{\rm
locu}(\overline{\D})$ and $g=f'$ are of the forms (\ref{structuref})
and (\ref{structureg}) and that $q(t)>0$ with $Q(t)\to\infty$ as
$t\to \infty$. Recall (\ref{rationalpoles}), (\ref{rationalpolesP})
and Corollary~\ref{cor:polesmoveout}, which in particular show that
the poles always move away from the origin. The following theorem
gives more precise estimates of their locations.

\begin{theorem}\label{first}
Assume $f(\zeta,t)$ is a global in time locally univalent solution
of the Polubarinova-Galin equation, with $g(\zeta,t)$ rational of
the form (\ref{structureg}). Then $|\zeta_j(t)|\sim a_1(t)$
for each $1\leq j\leq n$. More precisely, with $a_1(t)$ as in (\ref{powerseries})
$$
\frac{1}{a_1(0)}({|\zeta_j(0)|+\frac{1}{|\zeta_j(0)|}-2}) \leq
\frac{|\zeta_j(t)|}{a_1(t)} \leq
\frac{1}{a_1(0)}(|\zeta_j(0)|+\frac{1}{|\zeta_j(0)|}+2)
$$
for all $t\geq 0$.
\end{theorem}

\begin{proof}
Setting $\tau=\log a_1(t)$, $\xi=\xi(\tau)=|\zeta_j(t)|>1$, the
Harnack estimates (\ref{harnack}) say that
$$
\frac{\xi(\xi-1)}{\xi+1}\leq \frac{d\xi}{d\tau}\leq
\frac{\xi(\xi+1)}{\xi-1}.
$$
On integrating the differential equations corresponding to the equality cases one easily obtains the
inequalities
$$
\frac{1}{a_1(0)}({|\zeta_j(0)|+\frac{1}{|\zeta_j(0)|}-2})\leq \frac{1}{a_1(t)}({|\zeta_j(t)|+\frac{1}{|\zeta_j(t)|}-2})
$$
and
$$
\frac{1}{a_1(t)}({|\zeta_j(t)|+\frac{1}{|\zeta_j(t)|}+2})\leq \frac{1}{a_1(0)}({|\zeta_j(0)|+\frac{1}{|\zeta_j(0)|}+2}).
$$
Since
$$
\frac{1}{a_1(t)}({|\zeta_j(t)|+\frac{1}{|\zeta_j(t)|}-2})\leq
\frac{|\zeta_j(t)|}{a_1(t)}\leq \frac{1}{a_1(t)}({|\zeta_j(t)|+\frac{1}{|\zeta_j(t)|}+2})
$$
the desired estimates follow.

\end{proof}

Next, turning to coefficients we shall write $f$ and $g$ as follows:
\begin{equation}\label{structuref1}
f(\zeta,t)
=\sum_{j=1}^\ell e_j\log(1-\frac{\zeta}{\zeta_j(t)})+\frac{\sum_{j=1}^{m-\ell+1}b_j(t) \zeta^j}{\sum_{j=0}^{n-\ell} c_j(t)\zeta^j}
=\sum_{j=1}^\infty a_j(t)\zeta^j,
\end{equation}
\begin{equation}\label{structureg1}
g(\zeta,t)
=\frac{\sum_{j=0}^{m}\tilde{b}_j(t) \zeta^j}{\sum_{j=0}^{n} \tilde{c}_j(t)\zeta^j}
=\sum_{j=0}^\infty (j+1)a_{j+1}(t)\zeta^j.
\end{equation}
Introduce also the Taylor coefficients $\tilde{a}_j$ of the pure rational part of $f(\zeta,t)$ by
\begin{equation}\label{atilde}
\frac{\sum_{j=1}^{m-\ell+1}b_j(t) \zeta^j}{\sum_{j=0}^{n-\ell} c_j(t)\zeta^j}
=\sum_{j=1}^\infty \tilde{a}_j(t)\zeta^j.
\end{equation}

The above expansions are to be compared with (\ref{structureg}) and (\ref{structuref}).
In particular, $\ell$ denotes the number of different (finite) poles of $g$, $n_j$ denotes the order
of the pole at $\zeta_j$ (as in (\ref{structureg})) and
$$
n=\ell +\sum_{j=1}^\ell (n_j-1)=\sum_{j=1}^\ell n_j
$$
is the total order of the finite poles of $g$. The numbers $e_j$
are the residues of $g(\zeta)d\zeta$ at the points $\zeta_j$, and some or all of them may be zero.

At infinity,  $g(\zeta)d\zeta$ has a pole of order $n_0=m-n+2$ (as a differential), hence $f$ has a pole of order $m-n+1$ there.
This is also what (\ref{structuref1}) gives, hence the notations in (\ref{structuref1})
and (\ref{structureg1})
are consistent with those in (\ref{structuref}) and (\ref{structureg}).
Recall also that the coefficients $e_j$ do not depend on $t$. The
coefficients $b_j$, $c_j$, $\tilde{a}_j$, $\tilde{b}_j$, $\tilde{c}_j$ satisfy
$b_{m-\ell+1}\ne 0$, $c_{n-\ell}\ne 0$, $b_1\ne 0$, $c_0\ne 0$,
$\tilde{b}_{m}\ne 0$, $\tilde{c}_{n}\ne 0$, $\tilde{a_1}\ne 0$, $\tilde{b}_1\ne 0$, $\tilde{c}_0\ne 0$, and we shall
normalize them so that $c_0=1$, $\tilde{c}_0=1$.

\begin{lemma}\label{second}
Assume that $f(\zeta,t)$, given by (\ref{structuref1}), is a global
solution. Then, as $t\to \infty$,
$$
|c_{n-\ell}|\sim a_1^{-(n-\ell)}, \quad |c_j|= O(a_1^{-j}) \quad {\rm for\,\,}
1\leq j\leq n-\ell-1,
$$
$$
|b_1|\sim \tilde{a}_1\sim a_1, \quad |b_j|=O(1) \quad {\rm for\,\,} 2\leq j\leq m-\ell+1,
$$
$$
|\tilde{c}_n|\sim a_1^{-n}, \quad |\tilde{c}_j|= O(a_1^{-j}) \quad {\rm for\,\,}
1\leq j\leq n-1,
$$
$$
|\tilde{b}_0|\sim a_1, \quad |\tilde{b}_j|=O(1) \quad {\rm for\,\,} 1\leq j\leq m.
$$
\end{lemma}

See Subsection~\ref{sec:notations} for the meaning of $\sim$.

\begin{proof}

Since $\sum_{j=0}^{n-\ell} c_j(t)\zeta^j=c_{n-\ell}(t)\prod_{j=1}^{n-\ell}
(\zeta-\zeta_j(t))$ and $c_0=1$ the estimates for $c_j(t)$ follow
immediately from Theorem~\ref{first}. Similarly for $\tilde{c}_j(t)$.

Since by (\ref{structuref1}), (\ref{atilde}),
$$
{a}_j(t)=\tilde{a}_j(t)-\sum_{k=1}^\ell\frac{e_{k}}{j\zeta_k(t)^j},
$$
we have
$$
|\tilde{a}_j(t)-{a}_j(t)|\leq
\frac{C_j}{a_1(t)^{j}}\to 0,
$$
as $t\to\infty$. The coefficients $b_k$ are given by
\begin{equation}\label{bound0}
b_k=\tilde{a}_1 c_{k-1}
+\sum_{j=2}^{k-1} \tilde{a}_j c_{k-j}+\tilde{a}_k,
\end{equation}
and since the $a_j$, and hence the $\tilde{a}_j$, $j\geq 2$, are bounded (see (\ref{akN}) for example)
the assertions about the $b_j$ follow easily.
For the $\tilde{b}_k$ we have similarly
\begin{equation}\label{bound1}
\tilde{b}_k={a}_1 \tilde{c}_{k}
+\sum_{j=1}^{k-1} (j+1){a}_{j+1} \tilde{c}_{k-j}+(k+1){a}_{k+1},
\end{equation}
hence the estimates for these follow in the same way.

\end{proof}


\begin{theorem}\label{third}
Assume that $f(\zeta,t)$ in (\ref{structuref1}) is a global
solution and introduce the truncations
$$
f_{N}(\zeta,t)=\sum_{j=1}^{N}a_{j}(t)\zeta^{j}.
$$
Then the following assertions hold.
\begin{enumerate}

\item[(i)] There exist numbers $s_N$ with $s_{N}\rightarrow\infty$ as $N\rightarrow\infty$
such that, for $j=0,1$,
\begin{equation}\label{know}
\sup_{\zeta\in\overline{\mathbb{D}}}
\left|f_{N}^{(j)}(\zeta,t)-f^{(j)}(\zeta,t)\right|
=O(a_{1}(t)^{-s_{N}}),
\end{equation}
as $t\rightarrow\infty$. The same is true for any $j\geq 0$, with $s_N$ then depending on $j$.

\item[(ii)] For each $k\geq 2$,
\begin{equation}\label{decay1}
\lim_{t\rightarrow\infty}a_{k}(t)a_{1}(t)^{k}=\overline{M}_{k-1}.
\end{equation}

\item[(iii)] Assume that, for some number $r\geq 2$, $M_1=\dots =M_{r-1}=0$, $M_r\ne 0$. Then
\begin{equation}\label{decay2}
\lim_{t\rightarrow\infty}a_{s}(t)a_{1}(t)^{r+1}=0,\quad 2\leq s\leq r,
\end{equation}
and $r\leq m$. In case there are no logarithmic singularities (i.e., $e_j=0$, $1\leq j\leq \ell$) then
we even have $r\leq m-\ell$.

\end{enumerate}

\end{theorem}

\begin{proof}
Write $f$ on the form
$$
f(\zeta,t)=\sum_{j=1}^{\ell}e_{j}\log (1-\frac{\zeta}{\zeta_{j}})+b_1 \zeta\frac{1+P(\zeta)}{1-Q(\zeta)}
$$
$$
=-\sum_{k=1}^\infty\left(\sum_{j=1}^{\ell}\frac{e_j}{k\zeta_j^k}\right)\zeta^k
+b_1\zeta\sum_{k=0}^\infty(Q(\zeta)^k+P(\zeta)Q(\zeta)^k),
$$
namely with
$$
P(\zeta)=\sum_{j=1}^m \frac{b_{j+1}}{b_1}\zeta^j,
$$
$$
Q(\zeta)=-\sum_{j=1}^n c_j\zeta^j.
$$
Then all coefficients in $P$ and $Q$ are $O(a_1^{-1})$ by Lemma~\ref{second}, hence
$$
\sup_{\D(0,R)}|P|\leq \frac{C_R}{a_1}, \quad \sup_{\D(0,R)}|Q|\leq \frac{C_R}{a_1}
$$
for any fixed $R>1$ and suitable constants $C_R$.
It follows from the above that for large $N$ the remainder
$$
f(\zeta)-f_N(\zeta)=\sum_{j=N+1}^\infty a_j\zeta^j
$$
is built up by terms of the kind $\zeta Q(\zeta)^k$ and $\zeta P(\zeta)Q(\zeta)^k$
with also $k$ large. Considering only $t$ so large that $\frac{C_R}{a_1(t)} \leq \frac{1}{2}$
(for example) it follows that
\begin{equation}\label{supf}
\sup_{\D(0,R)}|f-f_N|\leq \frac{C_R}{a_1^{s}},
\end{equation}
where $s$ can be made arbitrarily large by choosing $N$ sufficiently large.
Now (\ref{know}) follows.

$(ii)$ Recall the expression (\ref{defmoments}) for the moments $M_k=M_k(f)$ in
terms of $f$. The truncations $f_N$ of $f$ similarly define moments:
\begin{align*}
M_{k}(f_N) &=\frac{1}{2\pi
i}\int_{\partial\D}f_{N}^{k}(\zeta,t)f_{N}'(\zeta,t)
{f_{N}^*}({\zeta},t)d\zeta\notag\\
&=\sum_{\sum i_j\leq N} i_{1}a_{i_{1}}\cdots
a_{i_{k+1}}\overline{a_{i_{1}+\cdots+i_{k+1}}}.
\end{align*}
It is clear from (\ref{defmoments}) that each moment $M_k(f)$ is a Lipschitz continuous function
of $f$ if $f$ is measured by the norm $\sup_{\overline{\D}}(|f|+|f'|)$, or any norm $\sup_{\D(0,R)}|f|$, $R>1$.
Therefore we have, by choosing $N=N_k$ large enough in (\ref{know}), with $j=0,1$
and on using (\ref{supf}), that for any exponent $p\geq 0$ we have
$$
|M_k(f)-M_k(f_N)|\leq C'|f-f_N|\leq \frac{C''}{a_1^p}
$$
for some constants $C'$ and $C''$ (which depend on $f$ and $p$).
From this we get (with $M_k=M_k(f)$)
$$
a_{1}^{k+1}\overline{a_{k+1}}-M_{k}
$$
$$
=M_{k}(f_N)-M_{k}-\sum_{k+2\leq \sum i_j\leq N}
i_{1}a_{i_{1}}\cdots a_{i_{k+1}}\overline{a_{i_{1}+\cdots+i_{k+1}}}
$$
\begin{equation}\label{estimate00}
= O(\frac{1}{a_{1}^p})-\sum_{k+2\leq\sum i_j\leq N} i_{1}a_{i_{1}}\cdots
a_{i_{k+1}}\overline{a_{i_{1}+\cdots+i_{k+1}}}
\end{equation}
$$
= O(\frac{1}{a_{1}^p})-(k+2)a_1^k a_2\overline{a_{k+2}}-\sum_{k+3\leq\sum i_j\leq N} i_{1}a_{i_{1}}\cdots
a_{i_{k+1}}\overline{a_{i_{1}+\cdots+i_{k+1}}}.
$$

Now, to prove (\ref{decay1}) we shall prove by
induction that for every $s\geq 2$, the two assertions
\begin{equation}\label{math01}
\lim_{t\rightarrow\infty}a_{1}(t)^{s}\overline{a_{s}(t)}=M_{s-1},
\end{equation}
\begin{equation}\label{math02}
|a_{j}(t)|=O(\frac{1}{a_{1}(t)^{s}}), \quad j\geq s
\end{equation}
hold. First, on using  the fact that the $|a_{j}(t)|$,
$j\geq 2$, are uniformly bounded we deduce from (\ref{estimate00}) that
\begin{equation*}
\left|a_{1}^{k+1}\overline{a_{k+1}}\right|\leq
|M_{k}|+O(\frac{1}{a_{1}^p}) +O(a_{1}^{k}), \quad k\geq 1,
\end{equation*}
and hence that $|a_{k+1}|=O(\frac{1}{a_{1}})$, $k\geq 1$.
Now repeating (\ref{estimate00}) with the new estimate
$|a_{j}| =O(\frac{1}{a_{1}})$, $j\geq 2$, we
obtain
\begin{equation*}
\left|a_{1}^{k+1}\overline{a_{k+1}}\right|\leq
|M_{k}|+O(\frac{1}{a_{1}^p}) +O(a_{1}^{k-2}), \quad k\geq 1,
\end{equation*}
and hence $|a_{k+1}|=O(\frac{1}{a_{1}^2})$, $k\geq 1$.

To start the induction process, take $k=1$ and use the last estimate
$|a_{j}|=O(\frac{1}{a_{1}^2})$, $j\geq 2$, in (\ref{estimate00}). This gives
\begin{equation*}
a_{1}^2\overline{a_{2}}-M_{1}=O(\frac{1}{a_{1}^p})+O(\frac{1}{a_{1}^3})
\end{equation*}
and hence $\lim_{t\rightarrow\infty}a_{1}^2\overline{a_{2}}=M_{1}$.
Therefore (\ref{math01}) and (\ref{math02}) hold for $s=2$.

Now take $s_0\geq 2$ and assume that (\ref{math01}), (\ref{math02}) hold for all $s\leq s_{0}$.
Then  we shall prove (\ref{math01}), (\ref{math02}) for $s=s_0+1$.
Thus we may in (\ref{estimate00}) use the fact that for any $s\leq s_0$ we have
$|a_{j}|=O(\frac{1}{a_{1}^{s}})$, $j\geq s$. This gives
\begin{equation*}
|a_{1}^{k+1}a_{k+1}|\leq
|M_{k}|+O(\frac{1}{a_{1}^p})+O(a_{1}^{k-2-s_{0}}), \quad k\geq s_{0},
\end{equation*}
and hence $|a_{k+1}|=O(\frac{1}{a_{1}^{s_{0}+1}})$ for $k+1\geq s_{0}+1$.
It follows that (\ref{math02}) holds for $s=s_{0}+1$.

Using, in
(\ref{estimate00}) with $k=s_{0}$, that
$|a_{j}|=O(\frac{1}{a_{1}^{s_{0}+1}})$ for $j\geq s_{0}+1$ (just proved) and
$|a_{j}|=O(\frac{1}{a_{1}^{j}})$ for $2\leq j\leq s_{0}$ (induction hypothesis), we
obtain
\begin{equation*}
a_{1}^{s_{0}+1}\overline{a_{s_{0}+1}}=M_{s_{0}}+O(\frac{1}{a_{1}^{p}})+O(\frac{1}{a_{1}^{3}})
\end{equation*}
and hence
$\lim_{t\rightarrow\infty}a_{1}^{s_{0}+1}\overline{a_{s_{0}+1}}=M_{s_{0}}$.
Therefore also (\ref{math01}) holds for $s=s_{0}+1$. Thus $(ii)$ in the theorem is proved.

$(iii)$ Assuming now $M_1=\dots=M_{r-1}=0$,
(\ref{estimate00}) with $1\leq k\leq r-1$ gives
\begin{equation}\label{ak1}
\overline{{a}_{k+1}}=O(\frac{1}{a_1^{p+k+1}})-\sum_{k+2\leq \sum i_j\leq N} i_{1}\frac{a_{i_{1}}\cdots
a_{i_{k+1}}}{a_1^{k+1}}\overline{a_{i_{1}+\cdots+i_{k+1}}}
\end{equation}
$$
=O(\frac{1}{a_1^{p+k+1}})-(k+2)\frac{a_2}{a_1}\cdot\overline{{a}_{k+2}}-\sum_{k+3\leq \sum i_j\leq N} i_{1}\frac{a_{i_{1}}\cdots
a_{i_{k+1}}}{a_1^{k+1}}\overline{a_{i_{1}+\cdots+i_{k+1}}}
$$
with $p$ arbitrarily large.
Using (\ref{decay1}) and  choosing $p\geq 4$, the above gives
$$
a_{k+1}=O(\frac{1}{a_1^{p+k+1}})+O(\frac{1}{a_1^{k+5}})=O(\frac{1}{a_1^{k+5}}).
$$

With $k=r-1$ this gives (\ref{decay2}) for $s=r$. It also follows
for $s=r-2$. Assume now, as an induction hypothesis, that (\ref{decay2}) holds for $s_0\leq s\leq r$.
If $s_0\geq 3$ we then let $k=s_0-2$ in (\ref{ak1}), and use (\ref{decay1}) and the induction hypothesis.
This gives (\ref{decay2}) for $s=s_0-1$. Hence (\ref{decay2}) is proved.

The identity (\ref{bound1})
holds for $1\leq k\leq m$, but also for $k>m$ with the convention that $b_k=0$ when $k>m$, and $\tilde{c}_j=0$ when $j>n$.
Choosing $k=m+1$ then gives
$$
(m+2)a_{m+2}=-\sum_{j=1}^m (j+1){a}_{j+1}\tilde{c}_{m+1-j}.
$$
If $r\geq m+1$ then (\ref{decay2}) and Lemma~\ref{second} show that the right member is
$o(a_1^{-(r+1)}\cdot O(a_1^{-1}))=o(a_1^{-(r+2)})$, while (\ref{decay1}) shows that
the left member is exactly $O(a_1^{-(r+2)})$, if $M_r\ne 0$. This contradiction shows that $r\leq m$.

Similarly, using (\ref{bound0}) one obtains $r\leq m-\ell$ in case there are no logarithmic terms
in $f$. An alternative way of proving these upper bounds for $r$ is given in the remark below.

\end{proof}

\begin{remark}
The fact that for  rational $g$, the vanishing of a sufficiently long sequence $M_1, M_2, \dots, M_k$ of moments
implies the vanishing of all $M_j$, $j\geq 1$, can also
be deduced from the relevant quadrature identity, like (\ref{qi}), holding for the image domain.
Assume for example that
$g(\zeta)d\zeta$ has no residues, so that there are no line integrals in (\ref{qi}). Then
choosing $h(\zeta)=\frac{1}{z-\zeta}$ with $z\in\C\setminus\overline{\Omega}$ in (\ref{qi}) gives the identity
$$
\sum_{k=0}^\infty\frac{M_k}{z^{k+1}}=
\sum_{j=0}^\ell \sum_{k=1}^{n_j-1}\frac{k!a_{jk}}{(z-z_j)^{k}}=\frac{\sum_{j=0}^{m-\ell}B_jz^j}{\sum_{j=0}^{m-\ell+1}C_jz^j}
$$
for suitable $B_j$, $C_j$. Recall that $\sum_{j=0}^\ell (n_j-1)=m-\ell+1$,
see at (\ref{qi}) and (\ref{structuref}). If now $M_1=\dots M_{r-1}=0$, $M_r\ne 0$, then
we get
$$
(M_0+\frac{M_r}{z^r}+\frac{M_{r+1}}{z^{r+1}}+\dots)(C_{m-\ell+1}+ \frac{C_{m-\ell}}{z} +\dots + \frac{C_0}{ z^{m-\ell+1}})
$$
$$
= B_{m-\ell}+\frac{B_{m-\ell-1}}{z}+\dots +\frac{B_0}{z^{m-\ell}}.
$$
Here it is easy to see from the general structure (\ref{structureg}), (\ref{structuref}) that
$M_0\ne 0$,  $C_{m-\ell+1}\ne 0$,  $B_{m-\ell}\ne 0$. On the other hand,
$C_{0}=C_1=\dots=C_{n_0-1}=0$ because the rational function above
has a pole of order (exactly) $n_0-1$ at $z_0=0$. Thus, in view of the
fact that $m-\ell +1-(n_0-1)=n-\ell$ the left
hand side of the above equation actually is
$$
(M_0+\frac{M_r}{z^r}+\frac{M_{r+1}}{z^{r+1}}+\dots)(C_{m-\ell+1}+ \frac{C_{m-\ell}}{z} +\dots + \frac{C_{n_0}}{ z^{m-\ell+1-(n_0-1)}})
$$
$$
= M_0 C_{m-\ell+1}+ \frac{M_0 C_{m-\ell}}{z} +\dots + \frac{M_0 C_{n_0}}{ z^{n-\ell}}+\frac{C_{m-\ell+1}}{z^r}+O(\frac{1}{z^{r+1}})+\dots.
$$
Comparing with the right hand and using that $n-\ell\leq m-\ell$ it follows that $r$ cannot be larger than
$m-\ell$, which is also what Theorem~\ref{third} tells.

In a similar way one proves that $r\leq m$ in the presence of logarithmic terms in $f$.

\end{remark}


From Theorem~\ref{third}, we obtain the following result:

\begin{cor}\label{four}
Let $f(\zeta,t)$ be a global solution of the form
(\ref{structuref1}) and assume
$M_1=\dots=M_{r-1}=0$, $M_r\ne 0$. Then,
\begin{equation*}\label{asym}
\lim_{t\rightarrow\infty}\sup_{\zeta\in\partial \mathbb{D}}
\left|\left[f(\zeta,t)-\sqrt{2Q(t)+M_{0}(0)}\zeta\right]
\left(\sqrt{2Q(t)}\right)^{r+1}-\overline{M_{r}}\zeta^{r+1}\right|=0.
\end{equation*}
\end{cor}

\begin{proof}
By (\ref{N0}), (\ref{diffu}), (\ref{diffusecond}) and Theorem~\ref{third},
\begin{equation*}\label{fin}
a_{1}(t)-\sqrt{2Q(t)+M_{0}(0)}=\frac{-N_0(t)}{a_{1}(t)+\sqrt{2Q(t)+M_{0}(0)}}
=O(\frac{1}{a_{1}^{2r+3}}).
\end{equation*}
Applying once more (\ref{diffu}), (\ref{diffusecond}) and Theorem~\ref{third}
the corollary follows.
\end{proof}


\bibliography{main0}

\end{document}